\documentclass{article}
\usepackage{amsmath}
\usepackage{amssymb}
\usepackage{amsfonts}
\numberwithin{equation}{section}
\usepackage[dvips]{graphicx}

\newtheorem{thm}{THEOREM}[section]
\newtheorem{prop}[thm]{PROPOSITION}
\newtheorem{lemma}[thm]{LEMMA}

\newtheorem{remark}[thm]{REMARK}
\numberwithin{figure}{section}

\newcommand{\dist}{\mathop{d_\C\,}}

\newcommand{\card}{\mathop{\rm card\,}}
\newcommand{\ind}{\mathop{\rm ind}}

\newcommand{\re}{\mathop{\rm Re\,}}
\newcommand{\e}{{\varepsilon}_0}
        
\newcommand{\calM}{{\mathcal M}}
\newcommand{\calN}{{\mathcal N}}
\newcommand{\calX}{{\mathcal X}}
\newcommand{\g}{{\gamma}}
\newcommand{\f}[1]{{(\ref{#1})}}
\newcommand{\rf}[1]{{\rm(\ref{#1})}}
\newcommand{\mea}{{\,dm_\D}}
\newcommand{\ber}{{B^2(\D,\bas)}}
\newcommand{\hil}{{{\mathfrak L}^2(\T,\bas)}}
\newcommand{\aff}{{L}}
\newcommand{\seg}{{\mathcal L}}
\newcommand{\Lop}{{\mathbf R}}
\newcommand{\calD}{{\mathcal D}}
\newcommand{\calR}{{\mathcal R}}
\newcommand{\calE}{{\mathcal E}}
\newcommand{\calO}{{\mathcal O}}
\newcommand{\calU}{{\mathcal U}}
\newcommand{\calQ}{{\mathcal Q}}
\newcommand{\Sop}{{\mathbf S}}
\newcommand{\Top}{{\mathbf T}}

\newcommand{\rdm}{{\theta}}
\newcommand{\bdm}{{\Theta}}
\newcommand{\rlm}{{\lambda}}
\newcommand{\blm}{{\Lambda}}
\newcommand{\bas}{{\omega}}
\newcommand{\bass}{{\widetilde{\omega}}}

\numberwithin{equation}{section}

\newcommand{\w}{{\Omega}}

\newcommand{\C}{{\mathbb C}}
\newcommand{\D}{{\mathbb D}}
\newcommand{\OD}{{\overline\D}}
\newcommand{\T}{{\mathbb T}}

\newcommand{\Z}{{\mathbb Z}}

\newcommand{\qed}{\vrule width 5pt height 5pt depth 0pt}

\newenvironment{proof}{\vspace{5pt}\bf Proof. \rm}{\hfill\qed\vspace{5pt}}
\newenvironment{proofof}[1]{\vspace{5pt}\bf Proof of #1. \rm}
{\hfill\qed\vspace{5pt}}

\begin{document}

\font\titlefont=cmr10 scaled \magstep3
\centerline{\titlefont{Large Bergman spaces:}}

\centerline{\titlefont{invertibility, cyclicity, and}}

\centerline{\titlefont{subspaces of arbitrary index}}

\bigskip
\bigskip
\bigskip

\centerline{\sl Alexander Borichev, H\aa{}kan Hedenmalm, and Alexander
Volberg}

\centerline{\sl at the University of Bordeaux I, Lund University,
and Michigan State University}

\bigskip
\bigskip
\bigskip
\bigskip



\begin{tabular}{p{12cm}} {\bf Abstract.}
In a wide class of weighted Bergman spaces, we construct invertible
non-cyclic elements. These are then used to produce $z$-invariant
subspaces of index higher than one. In addition, these elements generate 
nontrivial bilaterally invariant subspaces in anti-symmetrically weighted 
Hilbert spaces of sequences.
\end{tabular}
\bigskip
\bigskip
\bigskip

\section{\label{intro}Introduction}

Consider the operation which sends a
complex-valu\-ed sequence $\{a_n\}_n$ to the shifted sequence $\{a_{n-1}\}_n$;
if the index set is the collection of all integers, all is fine, and if it
is the nonnegative integers, we should specify that we need the rule
$a_{-1}=0$.
This shift operation is a linear transformation, and we denote it by $\Sop$.
Sometimes, it is convenient to work with formal Laurent or Taylor series
instead of sequence spaces, because of the simple form $\Sop$ takes, as it
just corresponds to the multiplication by the formal variable $z$:
$$\sum_n a_n z^n \,\,\mapsto \,\,\sum_n a_n z^{n+1}.$$
For an analyst, it makes sense to restrict the shift to some Hilbert space
like $\ell^2$, and then to use invariant subspaces to better understand
the operator. We recall that a linear subspace $\calM$ is {\sl shift
invariant} if it is closed with respect to the given topology, and if
$x\in\calM\Rightarrow \Sop x\in\calM$. The invariant subspaces offer
the possibility of studying the action of $\Sop$ on smaller pieces. However,
generally speaking, it may not
be possible to reconstruct the operator as a mosaic of really small and
understand\-able pieces. Arne Beurling \cite{Beu0} found a complete
characterization of the shift invariant subspaces in $\ell^2$ on the
nonnegative integers in terms of so-called inner functions. A consequence of
the theorem is that if a sequence (on the nonnegative integers) is in
$\ell^2$, and its convolution inverse is as well, then the smallest invariant
subspace containing it is all of $\ell^2$. Another consequence is that
every non-zero shift invariant subspace has index $1$; this means that
$\calM\ominus \Sop\calM$ is one-dimensional. If we instead
consider weighted $\ell^2$ spaces, with weights that make the space bigger
than unweighted $\ell^2$, we generally encounter a new kind of invariant
subspaces, having indices that exceed $1$; in fact, the dimension of
$\calM\ominus \Sop\calM$ can assume any (integer) value between $1$
and $+\infty$ \cite{ABFP}. Beurling later returned to shift invariant
subspaces in the setting of weights, but he softened the topology somewhat,
so instead of considering a weighted $\ell^2$ space, he looked at a union
of such, with the property of being an algebra, with respect to convolution
of sequences, which corresponds to ordinary multiplication of Taylor series. 
He found that there is a critical topology such that on the one side -- with
relatively small spaces -- the singular inner function for a point mass
generates a non-trivial invariant subspace, whereas on the other side of the 
borderline -- with relatively large spaces -- it generates the whole space
 as an invariant subspace. Later, Nikolai Nikolski \cite[Chapter 2]{Nik}
showed that the dichotomy is even deeper: we may replace the atomic singular
inner function by any zero-free function in the space. That is, in the
case of relatively large spaces, every zero-free function generates the whole 
space as an invariant subspace.

Here, we focus on hard topology spaces, that is, weighted $\ell^2$ spaces
on the nonnegative integers. We show that there is no such dichotomy as in
Beurling's situation. In fact, we find analytic functions in the disk that
belong to the given space along with their reciprocals, while the shift 
invariant subspaces they generate fail to be the whole space. We also use 
these functions to build concrete examples of invariant subspaces of high 
index.
Furthermore, we find that if we extend the weight and the space to the
collection of all integers, in such a way that the logarithm of the weight
becomes an odd function, then the nontrivial shift invariant subspace
generated by our invertible function extends to a bilaterally shift invariant 
subspace, in the sense that the intersection of the bilaterally invariant
subspace with the analytic part just returns us our initial shift invariant
subspace. 

\section{\label{descr-results}Description of results}

\subsection{\label{inv-cyc}Invertibility versus cyclicity}

In the Gelfand theory of commutative Banach algebras with unit,
an element generates a dense ideal if and only if it is
invertible, in which case its Gelfand transform has no zeros,
and the ideal it generates is the whole algebra. Let $\calX$ be a
Banach space (or a quasi-Banach space, that is, a complete normed space with a
$p$-homogeneous norm, for some $0<p<1$) of holomorphic
functions on the unit disk $\D$. We assume that the point
evaluations at points of $\D$ are continuous functionals on $\calX$,
and that $\Sop f\in \calX$ whenever $f\in \calX$, where $\Sop f(z)=z\,f(z)$ is
the operator of multiplication by $z$. This means that for each
$f\in \calX$, all the functions $\Sop f,\,\Sop^2f,\,\Sop^3f,\ldots$ are in
$\calX$ as well. Let $[f]=[f]_\calX$ denote the closure in $\calX$ of the
finite linear span of the vectors $f,\,\Sop f,\,\Sop^2f,\ldots$; we say
that $f$ is {\sl cyclic} in $\calX$ provided $[f]_\calX=\calX$. If
$\calX$ is a Banach (or quasi-Banach) algebra containing a unit element (the
constant function $1$), then all the polynomials belong to $\calX$.
If in fact the polynomials are dense in $\calX$, then $f\in \calX$ is
cyclic if and only if it is invertible. In other words, in {\sl
spaces} $\calX$, the concept of cyclicity generalizes that of
invertibility, provided that the polynomials belong to and are dense
in $\calX$. It is then of interest to compare cyclicity with genuine
invertibility. Consider, for instance, the space $\calX=H^2$, the
Hardy space on the unit disk $\D$. By Beurling's invariant
subspace theorem, a function is cyclic if and only if it is an
outer function. The invertible functions in $H^2$ are all outer,
so for this space, all invertible elements are cyclic. However,
by the examples provided by Borichev and Hedenmalm in
\cite{BH}, this fails for the Bergman spaces $B^p(\D)$,
$0<p<+\infty$, consisting of $p$-th power area-summable
holomorphic functions on $\D$. For an earlier example of a Banach space of
analytic functions where this phenomenon occurs, we refer to \cite{Sham}.

Given a continuous strictly positive area-summable function
$\bas$ on $\D$ --- referred to as a weight --- we form the space
$B^p(\D,\bas)$ (for $0<p<+\infty$) of holomorphic functions $f$
on $\D$ subject to the norm bound restriction
$$
\|f\|_{\bas,p}=\left(\int_\D|f(z)|^p\,\bas(z)\mea(z)\right)^{1/p}<
+\infty, 
$$ 
where $\mea(z)=\pi^{-1}\,dxdy$ is normalized area measure
($z=x+iy$). It is a Banach space of holomorphic functions for
$1\le p<+\infty$, and a quasi-Banach space for $0<p<1$. We shall
be concerned exclusively with radial weights $\bas$: from now
on, $\bas(z)=\bas(|z|)$ holds for all $z\in\D$. One nice thing
about radial weights is that the polynomials are guaranteed to
be dense in $B^p(\D,\bas)$. Let us form the soft topology space
${\mathcal A} (\D,\bas)=\cup_{0<p<+\infty}B^p(\D,\bas)$,
supplied with the inductive limit topology. Spaces of this kind were
studied by Beurling \cite{Beu} and Nikolski \cite[Chapter 2]{Nik}.
Under natural regularity conditions on $\bas$, the following dichotomy
holds: if
$$
\int_0^1\sqrt{\frac{\log\frac1{\bas(t)}}{1-t}}\,dt<+\infty, 
$$ 
then there exist a non-cyclic function in ${\mathcal A}
(\D,\bas)$ without zeros in $\D$ (the singular inner function
for an atomic measure will do), whereas if the above integral
diverges, that is,
$$
\int_0^1\sqrt{\frac{\log\frac1{\bas(t)}}{1-t}}\,dt=+\infty, 
$$ 
then each function in ${\mathcal A}(\D,\bas)$ lacking zeros in
$\D$ is cyclic in ${\mathcal A}(\D,\bas)$. One is then led to
wonder whether there exists a similar dichotomy for the hard
topology spaces $B^p(\D,\bas)$. Some progress has already been made on
this matter. Nikolski constructed in \cite[Section~2.8]{Nik} a special
class of weights $\bas$
that vanish at the boundary arbitrarily fast, and such that
$B^p(\D,\bas)$ contains zero-free noncyclic elements. Hedenmalm
and Volberg \cite{HV} proved that $B^2(\D,\omega)$, for
$$\omega(z)=\exp\left(-\frac1{1-|z|}\right),\qquad z\in\D,$$
contains invertible (and hence zero-free) noncyclic elements.
Aha\-ron Atzmon \cite{Atz} produced zero-free $S$-invariant subspaces
of index $1$ (for the notion of index see the next subsection) in $\ber$
for all $\bas$ satisfying some weak regularity conditions.

We formulate our first result. {\sl We shall assume that the
$($positive$)$ weight function $\bas$ decreases, and that $\omega(t)\to0$ as
$t\to1$ so quickly that for some $\e$, $0<\e<1$,}
\begin{equation}
\lim_{t\to 1}\,(1-t)^{\e}\,\log\log\frac1{\bas(t)}=+\infty;
\label{10}
\end{equation}
in other words, the speed is at least as fast as two exponentials. 
Without loss of generality, we can assume $\bas$ is $C^1$-smooth as well 
as
decreasing, and that the values are taken in the interval $(0,1/e)$.
No further assumptions of growth or regularity type will be made.

\begin{thm} For weights $\omega$ that meet condition {\rm(\ref{10})}, let
$\bass$ be the associated weight
\begin{equation}
\log\frac 1{\bass(z)}=\log\frac 1{\bas(z)}-
\Bigl[\log\log\frac 1{\bas(z)}\Bigr]^2,\qquad z\in\D,
\label{11}
\end{equation}
which decreases slightly more slowly than $\bas$ to $0$ as we
approach the boundary, so that $B^p(\D,\bass)$ is contained in
$B^p(\D,\bas)$. There exists a function $F\in B^1(\D,\bas)$
without zeros in $\D$ which is non-cyclic in $B^1(\D,\bas)$, and
whose reciprocal $1/F$ is in $B^1(\D,\bass)$. Moreover, we can
get $F$ such that in addition, $F^{1/p}$ is non-cyclic in
$B^p(\D,\bas)$ for each $p$, $0<p<+\infty$.
\label{t1}
\end{thm}

Note that both the growth and the decay of $|F|$ are somewhat
extremal. Indeed, were $F^{-1-\delta}$ to belong to $B^1(\D,\bas)$
for some $\delta>0$, then, for sufficiently regular $\bas$ we would
have $F^{1+\delta/2}\in B^1(\D,\bas)$. This follows from a uniqueness
theorem for harmonic functions \cite{Bor1}, which improves a result by
Nikolski \cite[Section~1.2]{Nik}). An elementary argument due to Harold 
Shapiro \cite{Sha} would then show that $F$ is cyclic in $B^1(\D,\bas)$.

It is interesting to note that we may strengthen the assertion of Theorem 
\ref{t1} to the following.

\begin{thm}
There exists a function $F$ satisfying the
conditions of Theorem~{\rm\ref{t1}} and such that for
$f=F^{1/2}$, we have
$$
\int_\D\int_\D \frac1{|f(z)|^2}\,\frac{|f(z)-f(w)|^2}{|z-w|^2}\,
\bas(z)\,\bas(w)\mea(z)\mea(w)<+\infty.
$$
\label{t1.5}
\end{thm}

\bigskip

\subsection{Subspaces of large index} 

The functions $F$ we construct are also extremal in a different
sense. Fix a $p$, $0<p<+\infty$. By the above theorem, there
exists a nontrivial invariant subspace of $B^p(\D,\bas)$ which
is generated by a function in $B^p(\D,\bas)$ whose reciprocal
belongs to the slightly smaller space $B^p(\D,\bass)$. Here we
use standard terminology: a closed linear subspace $\calM$ of
$B^p(\D,\bas)$ is {\sl invariant} if $\Sop f\in\calM$ whenever $f\in
\calM$. Note that for every $f\in B^p(\D,\bas)$, $\|\Sop f\|\ge C\,\|f\|$,
for some constant $C=C(p,\bas)$, $0<C<1$. Therefore, for every 
invariant subspace $\calM$,
the set $\Sop\calM$ is a closed subspace of $\calM$. The dimension
of the quotient space $\calM/\Sop\calM$ we shall call the {\sl index of}
$\calM$ and denote by $\ind \calM$. Given two invariant subspaces $\calM_1$
and $\calM_2$, we can form $\calM_1\vee \calM_2$, the smallest invariant
subspace containing both $\calM_1$ and $\calM_2$. This definition
naturally extends to more general collections of invariant
subspaces, with more than two elements. The index function is
then subadditive, in the following sense: the index of
$\calM_1\vee\ldots\vee \calM_n$ is less than or equal to the sum of the
individual indices for $\calM_1,\ldots, \calM_n$.

Constantin Apostol, Hari Bercovici, Ciprian Foia{\c s}, and Carl
Pearcy \cite{ABFP} proved that $\ber$ contains invariant
subspaces of arbitrary index. Theirs is a pure existence theorem,
but later on, concrete examples of invariant subspaces of large
index were given in \cite{Hed, HRS, Bor, AB}. Usually, subspaces
of index bigger than $1$ have somewhat strange properties. For
example, they cannot contain multipliers \cite{R}; in
$B^p(\D)$, $1<p<+\infty$, if $\ind[f,g]=2$ (where $[f,g]=
[f]\vee[g]$ is the closed invariant subspace generated by $f$
and $g$), then $f/g$ has finite non-tangential limits almost
nowhere on the unit circle \cite[Corollary~7.7]{ARS}. It turns
out that the invertible functions $F$ we construct could be used to 
produce invariant subspaces of higher index.

\begin{thm} Let us assume that the weight $\bas$ satisfies property 
\rf{10}. Then there exist invertible functions $F_j$ in $B^1(\D,\bas)$
for $j=1,2,3,\ldots$, such that for every $0<p<+\infty$ and for every
positive integer $n$, the subspace
$$
[F_1^{1/p},\ldots,F_n^{1/p}]=[F_1^{1/p}]\vee\ldots\vee[F_n^{1/p}]
$$
has maximal index in $B^p(\D,\bas)$, namely $n$. Moreover, the
assertion holds also when $n$ assumes the value $+\infty${\rm:}
$$
[F_1^{1/p},F_2^{1/p},\ldots]
$$
has infinite index in $B^p(\D,\bas)$.
\label{t2}
\end{thm}
\bigskip

\subsection{The induced bilateral Hilbert space}

The space $\ber$ is a Hil\-bert space, and it is possible to
describe it as a weighted $\ell^2$ space on the set of
nonnegative integers $\Z_+=\{0,1,2,\ldots\}$. At times, we also need
the collection of negative integers $\Z_-=\{-1,-2,-3,\ldots\}$.
We note that a holomorphic function
$$
f(z)=\sum_{n=0}^{+\infty} \widehat f(n)\,z^n,\qquad z\in\D,
$$
is in $\ber$ if and only if
\begin{equation}
\sum_{n=0}^{+\infty}|\widehat f(n)|^2\,\w(n)<+\infty,
\label{10.5}
\end{equation}
where
\begin{equation}
\w(n)=\int_\D |z|^{2n}\bas(|z|)\mea(z)<+\infty,\qquad n=0,1,2,3,\ldots.
\label{10a}
\end{equation}
The function $\w$ is log-convex, that is, $\w(n)^2\le
\w(n-1)\w(n+1)$ holds for all $n=1,2,3,\ldots$, and it has the
property
\begin{equation}
\lim_{n\to+\infty}\w(n)^{1/n}=1.
\label{10m}
\end{equation}

The left hand side of \f{10.5} equals the norm squared of $f$ in
$\ber$. With the usual Cauchy duality 
$$
\langle f,g\rangle =\sum_{n=0}^{+\infty} a_nb_n, 
$$ 
where $g$ is the convergent Laurent series 
$$
g(z)=\sum_{n=0}^{+\infty} b_n\, z^{-n-1},\qquad 1<|z|<+\infty, 
$$
we can identify the dual space $\ber^*$ with the space of
Laurent series $g$ with norm
\begin{equation}
\|g\|_{\bas*}^2= \sum_{n=0}^{+\infty}\frac{|b_n|^2}{\w(n)}<+\infty. 
\label{dualbergspace}
\end{equation}
We then form the sum space $\hil=\ber\oplus\ber^*$ of formal Laurent series 
$$
h(z)=\sum_{n=-\infty}^{+\infty} c_n\,z^n,\qquad |z|=1, 
$$
with norm
$$
\|h\|_{\hil}^2=\sum_{n=0}^{+\infty} |c_n|^2\w(n)
+\sum_{n=0}^{+\infty} \frac{|c_{-n-1}|^2}{\w(n)}. 
$$ 
Note that although we have indicated the unit circle 
as the domain of definition, in general the formal series 
converges nowhere. If we extend the weight $\w$ to negative integers by

$$
\w(n)=\frac1{\w(-n-1)},\qquad n=-1,-2,-3,\ldots,
$$
then the norm on $\hil$ takes on a more pleasant appearance:
$$
\|h\|_{\hil}^2=\sum_{n=-\infty}^{+\infty} |c_n|^2\,\w(n). 
$$ 

On the other hand, given a positive log-convex function $\w_0$
such that
$$\lim_{n\to+\infty}\w_0(n)^{1/n}=1,$$
we can find a weight function $\bas$ on $[0,1)$ such that the function $\w$
defined by $\bas$ via formula \f{10a} is equivalent to $\w_0$, written
out
$\w_0\asymp\w$ (see \cite[Proposition~B.1]{BHA}). In concrete terms, 
this
means that for some positive constant $C$,  $C^{-1}\w_0\le\w\le C\,\w_0$ 
on the nonnegative integers.
Let us see what the condition \f{10} requires in terms of the weight 
$\w$. Suppose we know that for some $0<\alpha<+\infty$,
\begin{equation}
\w(n)\le\exp\left[-\frac{n}{(\log(2+n))^\alpha}\right],
\qquad n=0,1,2,\ldots.
\label{10k}
\end{equation}
Then the decreasing weight $\bas$ automatically satisfies \f{10}. 
Indeed, by \f{10a} we have
$$
2\bas(x)\int_0^x r^{2n+1}\,dr\le
2\int_0^1 r^{2n+1}\bas(r)\,dr\le
\exp\left[-\frac{n}{(\log(n+2))^\alpha}\right],
$$
for $n=0,1,2,3,\ldots$, and hence
\begin{equation*}
\bas(x)\le \min_n\frac{(n+1)\exp\Big[-\frac{n}{(\log(n+2))^\alpha}\Big]}
{x^{2n+2}},\qquad 0<x<1,
\end{equation*}
from which \f{10} follows.

The shift operator $\Sop$ extends to $\hil$,
$$ 
\Sop h(z)=z\,h(z)=\sum_{n=-\infty}^{+\infty} c_{n-1}\,z^n, 
$$ 
but now $\Sop$ is an {\sl invertible} bounded operator. It is therefore
natural to consider closed subspaces that are invariant with respect
to both the forward shift and the backward shift $\Sop^{-1}$. We call
them bilaterally invariant, and the shift in both direction the
bilateral shift.

Let $\calN$ be a bilaterally invariant subspace of $\hil$, and $\calM$ a
subset of $\ber$. We consider the intersection $\calN_\D=\calN\cap
\ber$, and the extension $\calM_\T$ of $\calM$, the closed linear span
of $\bigcup_{n\in\Z}\Sop^n\calM$ in $\hil$. Then $\calN_\D$ is a forward
invariant subspace of index $1$ in $\ber$, and $\calM_\T$ is a
bilaterally invariant subspace of $\hil$.

If the weight $\bas$ is sufficiently regular and 
$$
\int^1\log\log\frac1{\bas(x)}\,dx<+\infty,
$$
or, equivalently,
$$
\sum_{n=0}^{+\infty}\frac{\log\frac1{\w(n)}}{n^2+1}<+\infty,
$$
then the classical result of John Wermer \cite{W} implies that
$\hil$ possesses nontrivial (bilaterally) invariant subspaces.
It seems that the existence of nontrivial bilaterally invariant
subspaces $\calN$ in $\hil$ is unknown for general $\bas$ (for the current
status, see \cite{AS,Atz2});
however, for sufficiently regular weights $\bas$, this was proved
by Yngve Domar \cite{Dom}.

In his example, we have $\calN_\D=\{0\}$. 
Let us sketch an argument to this effect.
Recall that $\hil$ is isometrically isomorphic to the space
$\ell^2(\Z,\w)$ of complex-valued sequences $\{c_n\}_{n\in\Z}$ with
$$
\big\|\{c_n\}_{n}\big\|_{\ell^2(\Z,\w)}^2=
\sum_{n\in\Z}|c_n|^2\w(n)<+\infty,
$$
provided that $\w$ is extended ``anti-symmetrically'' to the
negative integers:
$$
\w(n)=\frac1{\w(-n-1)},\qquad n<0.
$$
Note that
$$0<\inf_{n\in\Z_+}\frac{\w(n-1)}{\w(n)}
\le\sup_{n\in\Z_+}\frac{\w(n-1)}{\w(n)}<+\infty,$$
so that for all essential purposes, we can think of $\log\w(n)$ as an
anti-symmetric function of $n$.
Domar constructs a nontrivial entire function $f$ of exponential type
$a$, with $0<a<\frac12\pi$, such that
$$|f(n)|^2\le \frac1{(n^2+1)\w(n)},\qquad n\in\Z.$$
Then
$$
\widehat f(z)=\sum_{n\in\Z}f(n)\, z^n\in\hil,
$$ 
and $[\widehat f]_\T\ne\hil$ (here, of course, $[\widehat f]_\T$ is the
bilaterally invariant subspace generated by $\widehat f$ in $\hil$). To
show that $[\widehat f]_\T\cap\ber=\{0\}$, we 
suppose that for a sequence of entire functions $f_k$ of exponential
type $a$ we have $f_k{\bigm|}_\Z\to h$ in $\ell^2(\Z,\w)$, where
$h{\bigm|}_{\Z_-}\equiv 0$, $h(0)\ne 0$. Then $F_k(z)=f_k(z)f_k(-z)
\to H(z)=h(z)h(-z)$ in $\ell^1(\Z)$. Applying the Cartwright theorem
\cite[Section~21.2]{L} and the Phragm\'en--Lindel\"of theorem in
the upper and the lower half-planes, we obtain that as $k\to+\infty$, $F_k$
converges to an entire function $F$ of exponential type at most
$2a$, and that $F{\bigm|}_\Z=H{\bigm|}_\Z$. Note that $H(n)=0$ on
$\Z\setminus\{0\}$; but this requires $H$ to have at least type $\pi$, which
is not possible, by the assumption that $a<\frac12\pi$.
\medskip

In the related work \cite{EV}--\cite{EV2}, Jean Esterle and
Alexander Volberg proved that for asymmetrically weighted Hilbert
spaces of sequences $\ell^2(\Z,\w)$, with $\w(n)\ll \w(-n)^{-1}$ as
$n\to+\infty$, and
$$\sum_{n=1}^{+\infty}\frac{\log\w(-n)}{n^2}=+\infty,$$

$(i)$ every bilaterally invariant subspace $\calN$ is generated by its
analytic part $\calN\cap\ell^2(\Z_+,\w)$, where $\ell^2(\Z_+,\w)$ is the
subspace of $\ell^2(\Z,\w)$ consisting of all sequences vanishing on the
negative integers $\Z_-$, and

$(ii)$ $\calN+\ell^2(\Z_+,\w)=\ell^2(\Z,\w)$.
\medskip

On the other hand, for anti-symmetrically weighted $\ell^2(\Z,\w)$, with
sufficiently regular $\w$, such that for some positive value of the
parameter $\epsilon$,
$$n^{1+\epsilon}\le\w(n)$$
for all big $n$, Borichev \cite{Bnew} proved that for every bilaterally
invariant subspace $\calN$, $\calN\cap\ell^2(\Z_+,\w)=\{0\}$. 

Let us give here an application of our Theorem~\ref{t1}, which actually
requires a slightly stronger property of the function $F$, . If $\w$ is a
log-convex weight function satisfying the conditions of
Theorem~3-3 in \cite{EV1}, such that \f{10k} holds for some
$\alpha>0$, then $\ell^2(\Z_+,\w)$ contains a non-trivial zero-free
shift invariant subspace of index $1$; thus, by Theorem~3-3 of
\cite{EV1}, $\ell^2(\Z,\w)$ contains non-trivial bilaterally
invariant subspaces.

We construct bilaterally invariant subspaces in $\hil$  with
properties similar to $(i)$ and $(ii)$.

\begin{thm} Let us assume that the weight $\bas$ satisfies property 
\rf{10}. Then there exists a function $F$ satisfying the
conditions of Theorem {\rm\ref{t1}} and such that $F^{1/2}$ is
non-cyclic in $\hil$ with respect to the the bilateral shift.
Furthermore, if $\calM=[F^{1/2}]$, then
\begin{equation}
\left.
\begin{gathered}
(\calM_\T)_\D=\calM,\\
\calM_\T+\ber=\hil.
\end{gathered}
\right\}
\label{10b}
\end{equation}
\label{t3}
\end{thm}

\begin{remark} {\rm (a) In view of the theorem and the above observations,
if $\w$ is a log-convex weight function satisfying \f{10m} and \f{10k}
for some $\alpha>0$, and extended to negative indices $n$ by
$\w(n)=1/\w(-n-1)$, then $\ell^2(\Z,\w)$ contains a singly generated
bilaterally invariant subspace $\calN$ such that

$(i)$ $\calN$ is generated by $\calN\cap\ell^2(\Z_+,\w)$, and

$(ii)$ $\calN+\ell^2(\Z_+,\w)=\ell^2(\Z,\w)$.
\medskip

(b) It is interesting to contrast the situation depicted in part (a) with
the case $\w(n)\equiv1$, when each bilaterally shift invariant subspace
$\calN$ of $\ell^2(\Z)$ is given by a common zero set on the unit circle
$\T$, and we have $\calN\cap\ell^2(\Z_+)=\{0\}$ unless $\calN=\ell^2(\Z)$.
An intermediate case between these opposites would be the Bergman-Dirichlet
situation, with $\w(n)=1/(n+1)$ for $n=0,1,2,\ldots$ and $\w(n)=-n$ for
$n=-1,-2,-3,\ldots$. Here, it is not known whether there exist nontrivial
bilaterally invariant subspaces $\calN$ with property $(i)$ above.
}
\end{remark}
\bigskip

\subsection{The spectra associated with bilaterally invariant 
subspaces}

Given an invariant subspace $\calM$ in $\ber$, the operator $\Sop$
induces an operator $\Sop[\calM]:\ber/\calM\to\ber/\calM$ defined by
$\Sop[\calM](f+\calM)=\Sop f+\calM$; clearly, $\Sop[\calM]$ has norm less than
or equal to
that of $\Sop$, which equals $1$. Similarly, given a bilaterally
invariant subspace $\calN$ in $\hil$, the operator $\Sop$ induces an
operator $\Sop[\calN]:\hil/\calN\to\hil/\calN$ defined by
$\Sop[\calN](f+\calN)=\Sop f+\calN$.
Again, $\Sop[\calN]$ has norm less than or equal to that of $\Sop$, and
$\Sop[\calN]^{-1}$ has norm less than or equal to that of $\Sop^{-1}$. We
do not indicate the underlying space here in the notation,
because we feel that no confusion is possible. In the situation
indicated in Theorem \ref{t3}, the operators $\Sop[\calM]$ and
$\Sop[\calM_\T]$ are canonically similar. For, it is easy to check that
the canonical mapping
$$
j_\calM:\ber/\calM\to\hil/\calM_\T
$$ 
given by $j_\calM(f+\calM)=f+\calM_\T$ is an isomorphism, and we have the
relationship $\Sop[\calM_\T]=j_\calM\circ \Sop[\calM]\circ j_\calM^{-1}$. In
particular, the operators $\Sop[\calM]$ and $\Sop[\calM_\T]$ have the same
spectrum. The spectrum of $\Sop[\calM_\T]$ is a compact subset of $\T$,
because the unit circle is the spectrum of the bilateral shift
on $\hil$. Generally speaking, the spectrum of $\Sop[\calM]$ is the
common zero set of $\calM$ on the open unit disk plus a generalized
zero set on the unit circle (see \cite{Hed2}, which treats the
unweighted Bergman space case). It is a consequence of the next
theorem that the above situation depicted in Theorem \ref{t3}
cannot occur if $\Sop[\calM]$ has countable spectrum (contained in the
unit circle).

\begin{thm} Suppose $\calN$ is a bilaterally shift invariant
subspace of $\hil$. Then the spectrum $\sigma(\Sop[\calN])$ of $\Sop[\calN]$
is a perfect set, that is, a closed subset of $\T$ without
isolated points. In particular, if $\sigma(\Sop[\calN])$ is countable,
then $\sigma(\Sop[\calN])=\emptyset$, and $\calN=\hil$.
\label{t4}
\end{thm}

The proof of the above theorem in Section~\ref{pf-t4} does not use any of
the strong assumptions made earlier on the weight, but rather
holds in a much more general context. In view of the above theorem, it is
natural to ask what kinds of sets may actually occur as spectra of 
$\Sop[\calN]$, if $\calN$ is a bilaterally shift invariant subspace of $\hil$.
It is possible to verify that the bilaterally invariant subspace $\calM_\T$ 
appearing in Theorem \ref{t3} has spectrum $\T$, that is, the induced operator 
$\Sop[\calM_\T]$ has spectrum $\T$. But perhaps there are other bilaterally
invariant subspaces with more complicated spectra? At the present moment in
our investigation, we do not even know if the space $\hil$ always possesses
a bilaterally invariant subspace whose spectrum is a non-trivial closed arc
of $\T$. Possible candidates might be the subspaces constructed by Domar
\cite{Dom}. Once this question is resolved, it is natural
to ask what it means for the spectrum of a bilaterally invariant subspace 
$\calN$ if we add the requirement that the analytic part of $\calN$ is 
nontrivial, that is, $\calN\cap B^2(\D,\bas)\neq\{0\}$.

\bigskip

\subsection{The idea of the proof of Theorem~{\ref{t1}}}

The largeness of a function sometimes implies its non-cyclicity. At
first glance this is counter-intuitive, because in Hardy spaces,
for example, outer functions are the ``largest'' ones, but they
are cyclic. However, in spaces of analytic functions determined
in terms of growth, the largeness does imply non-cyclicity. One
can get a flavor of the proof by first looking at the following
``toy picture''. Let a class of analytic functions $f$ be
defined by the condition $|f(z)|\le C\,\varrho(|z|)$, where
$\varrho$ is a (radial) weight of growth type, that is
$\varrho(x){\nearrow}+\infty$, $x\to 1$.  Furthermore, let $F$ be
in this class, and let it be ``maximally large'' in the sense
that $|F(z)|\ge\varrho(|z|)$ for all $z\in\calQ\subset\D$.
If the set $\calQ$ is massive enough, here is what will
happen. Any sequence of polynomials $q_n$ such that
$|q_n(z)\,F(z)|\le C\,\varrho(|z|)$ will obviously satisfy the
estimate $|q_n(z)|\leq C$, $z\in\calQ$. Now the massiveness of
$\calQ$ guarantees the uniform boundedness of the family
$\{q_n\}_n$ (here massiveness may mean, for example, that almost
every point of the circle can be approached by points from
$\calQ$ in a non-tangential way). The uniform boundedness
of the family $\{q_n\}_n$ should be combined with a property of
$F$ (which one has to establish in advance) to tend to zero
along a certain sequence. Together, these two properties show
that the products $q_n\,F$ cannot converge to a nonzero constant
uniformly on compact subsets of the disk. Thus, the non-cyclicity
of $F$ follows.

This kind of idea was used to construct non-cyclic functions in
the paper of Borichev and Hedenmalm \cite{BH}. The same idea
will also be used in the present article. However, we shall not
be able to prove the uniform boundedness of $q_n$. Instead, we
shall prove the normality of the family $\{q_n\}_n$, with
effective uniform estimates on the growth of $|q_n(z)|$.
The difference with the ``toy problem'' is that now, $F$ will be
of maximal largeness in the integral sense, that is, on average,
rather than pointwise as above. In its turn, this will imply
that $q_n$ are not uniformly bounded on a massive set as before,
but rather have some integral estimates; more precisely, the weighted
sum of the absolute values of $q_n$ along a lattice-like sequence of points
in $\D$ will satisfy some effective estimates (independent of $n$).
Unlike in the previous consideration, where the
massiveness of $\calQ$ was used (essentially) via the
harmonic measure estimates, we will use the Lagrange
interpolation theorem.

To describe the idea of ``integral largeness'' in more details,
we start with the weight $\omega$ of decrease type from Subsection
\ref{inv-cyc}.
We will construct $F$ in the unit ball of $B^1(\mathbb{D},
\omega)$, which is maximally large in a certain sense.
What this amounts to is
$$
\int_\D |F(z)|\,\omega(|z|)\mea(z)\leq 1
$$
and (the inequality below is what we mean by ``integral maximal
largeness'')
$$
\int_{\calD_{n,k}}|F(z)|\,\omega(|z|)\mea(z)
\asymp\frac {1-r_n}{n^2}, 
$$
where $r_n$ is a sequence of radii rapidly converging to $1$,
and $w_{n,k}$, with $0\le k<N_n\asymp\frac1{1-r_n}$,
are the points equidistributed over the circle of radius $r_n$ about the
origin, and finally, $\calD_{n,k}$ are the disks
$$\calD_{n,k}=\{z\in\C:\,|z-w_{n,k}|<(1-r_n)^2\}.$$
Summing up over $k,n$, we notice that the integral
of $|F(z)|\,\omega(|z|)$ over the union of our small disks
is proportional to the integral over the whole unit disk. A big
part of the mass of $|F(z)|\,\omega(|z|)\mea(z)$ lies
inside a tiny set, which is the union of small disks. In other
words, a big part of the norm of $F$ is concentrated on a set
which is tiny in the sense of area, but which is sufficiently
``widespread''.

Suppose that for a sequence of polynomials $q_m$,
$$
\|q_m\,F\|_{B^1(\bas)}\le A,
$$
where $A$ is a positive constant. Then the ``integral maximal largeness''
inequality shows that for some points $z_{n,k}$ (the points will depend on
the choice of $q_m$, they are the points of minimum for $|q_m(z)|$ on the
corresponding small disks $\calD_{n,k}$) we have the weak type estimate
$$
\sum_{k=0}^{N_n-1}|q_m(z_{n,k})|\le C\,n^2N_n.
$$
Applying the results of Section~\ref{main-constr}, we get for some
positive constant $C$ independent of $m$ that
$$
|q_m(z)|\le C\,\exp\left[\frac1{1-|z|}\right],\qquad z\in\D.
$$
This is the effective estimate on normality we mentioned
above. Next, we have to guarantee that $F$ tends to zero
sufficiently rapidly along a sequence of points. This, together
with the last estimate, implies that the products $q_m \,F$
cannot converge to a nonzero constant uniformly on compact
subsets of the disk. Thus, the non-cyclicity of $F$ follows.

\bigskip

\subsection{The plan of the paper}

First, in Section~\ref{pf-t4}, we prove Theorem~\ref{t4}. Then, in Section
\ref{pf-t3}, we show how to use Theorem \ref{t1.5} to deduce Theorem \ref{t3},
applying a metod similar to that used in \cite{EV1}. After this, we turn to 
the more technical aspects of the paper. In Section~\ref{main-constr}, we 
construct some elements in $B^1(\D,\bas)$ of ``maximal possible growth''. 
The constructions use 

\noindent $(a)$ a regularization procedure for $\bas$ described in 
Section~\ref{aux-convex},

\noindent $(b)$ estimates for auxiliary harmonic functions established in
Sections~\ref{aux-harm-constr} and \ref{aux-harm-est}, and 

\noindent $(c)$ a special Phragm\'en--Lindel\"of type estimate obtained in 
Section~\ref{phrlind-est}.

Theorems~\ref{t1}, \ref{t1.5}, and \ref{t2} are proved in 
Section~\ref{main-constr}, renamed as Theorems~\ref{p1}, \ref{p2}, and 
\ref{p3}, respectively.
\bigskip
\bigskip

\section{\label{pf-t4}The proof of Theorem \ref{t4}}

The spectrum $\sigma(\Sop[\calN])$ is a closed subset of the unit
circle $\T$. The set $\sigma(\Sop[\calN])$ is empty if and only if
$\Sop[\calN]$ acts on the trivial space $\{0\}$, in which
case $\calN=\hil$. A nonempty closed and countable subset of $\T$
necessarily has isolated points. It remains to prove that
isolated points cannot occur in the spectrum $\sigma(\Sop[\calN])$. To
this end, let $\lambda_0\in\sigma(\Sop[\calN])$ be an isolated point.
By the Riesz decomposition theorem (see, for example, 
\cite[Section~2.2]{RR}), the bilaterally invariant subspace
$\calN$ can be written as $\calN=\calN_0\cap \calN_1$, where
$\sigma(\Sop[\calN_0])=
\{\lambda_0\}$ and $\sigma(\Sop[\calN_1])=\sigma(\Sop[\calN])\setminus
\{\lambda_0\}$. We shall prove that $\Sop[\calN_0]$ cannot have
one-point spectrum, which does it. In other words, after
replacing $\calN$ by $\calN_0$ and after a rotation of the circle, we
may assume that $\calN$ has $\sigma(\Sop[\calN])=\{1\}$. We use the
holomorphic functional calculus (again, in fact) to define the
operator $\log \Sop[\calN]$, with spectrum $\{0\}$. This permits us to
form the expression
$$
\Sop[\calN]^\zeta=\exp\big(\zeta\log \Sop[\calN]\big),\qquad \zeta\in\C,
$$
which amounts to an entire function of zero exponential type
taking values in the space of operators on $\hil/\calN$. We identify
the quotient space $\hil/\calN$ with the subspace $\hil\ominus \calN$ of
vectors perpendicular to $\calN$. Let ${\mathbf P}:\hil\to\hil
\ominus \calN$ stand for the orthogonal projection, and denote by
$$
\langle f,g\rangle=\sum_{n\in\Z}\widehat
f(n)\,\overline{\widehat g(n)}\,\w(n)
$$
the sesquilinear form on the Hilbert space $\hil$, where for a
given element $f\in\hil$, $\widehat f(n)$ are the corresponding
Laurent coefficients for the formal series expansion
$$
f(z)=\sum_{n\in\Z}\widehat f(n)\,z^n,\qquad z\in\T.
$$
We recall the agreed convention that $\w(n)=1/\w(-n-1)$ for $n<0$.
The operator $\Sop[\calN]$ is identified with ${\mathbf P}\Sop$. For
$f,\phi\in\hil$, we consider the function
$$
E_{f,\phi}(\zeta)=\langle \Sop[\calN]^\zeta\,{\mathbf P}f,
{\mathbf P}\phi\rangle,\qquad \zeta\in\C,
$$
which is entire and of zero exponential type. Using the
invariance of $\calN$, we see that ${\mathbf P}\Sop{\mathbf P}={\mathbf
P}\Sop$; more generally, for integers $n$, we have
$$
\Sop[\calN]^n{\mathbf P}=({\mathbf P}\Sop)^n{\mathbf P}={\mathbf P}\Sop^n.
$$
holds, so that as we plug in $f=1$ and $\phi\in\hil\ominus \calN$
into the above expression, we obtain, for $n\in\Z$,
$$
E_{1,\phi}(n)=\langle \Sop[\calN]^n\,{\mathbf P}1,
\phi\rangle=\langle{\mathbf P}\Sop^n 1,\phi\rangle=
\langle z^n,\phi\rangle=\w(n)\,\overline{\widehat\phi(n)}.
$$
On the other hand, if we instead fix $\phi=1$ and let
$f\in\hil\ominus \calN$ vary, we obtain
$$
E_{f,1}(\zeta)=\langle \Sop[\calN]^\zeta\,f,{\mathbf P}1\rangle=
\langle \Sop[\calN]^\zeta\,f,1\rangle,
$$
and hence
$$
E_{f,1}(n)=\langle z^n\,f,1\rangle=\w(0)\,\widehat f(-n),
\qquad n\in\Z.
$$
Now, the entire function $F(\zeta)=E_{f,1}(-\zeta)
E_{1,f}(\zeta)$ has zero exponential type, and by the above, it
is $l^1$-summable at the integers:
\begin{multline*}
\sum_{n\in\Z}|F(n)|=
\sum_{n\in\Z}\big|E_{f,1}(-n)E_{1,f}(n)\big|\\
=\w(0)\sum_{n\in\Z}|\widehat f(n)|^2\w(n)
=\w(0)\|f\|_{\hil}^2<+\infty.
\end{multline*}
By the classical Cartwright theorem, the function $F$ is bounded
along the real line, and in view of the growth restriction which
is a consequence of $F$ having zero exponential type, the
Phragm\'en--Lindel\"of principle forces $F$ to be constant. That
constant must then be $0$, in view of the convergence of the
above sum. Then at least one of the entire functions $E_{f,1}$
and $E_{1,f}$ must vanish identically. In either case, $\widehat
f(n)=0$ for all integers $n$, that is, $f=0$. Since $f$ was
arbitrary in $\hil\ominus \calN$, we obtain $\calN=\hil$, and hence
$\sigma(\Sop[\calN])=\emptyset$, as desired.
\qed

\bigskip
\bigskip

\section{\label{pf-t3}The proof of Theorem \ref{t3}}

We start with a function $f\in\ber$ as in Theorem~\rm\ref{t1.5}.
So, we know that $[f]\ne\ber$ and
\begin{equation}
\int_\D\int_\D \frac1{|f(z)|^2}\frac{|f(z)-f(w)|^2}{|z-w|^2}
\,\bas(z)\,\bas(w)\mea(z)\mea(w)<+\infty.
\label{44a}
\end{equation}
Recall that $\Sop$ is the shift operator on $\hil$, and
consider the induced operator $\Top$ on $\ber/[f]$. Since
$\Sop$ is a contraction on $\ber$, we have $\|\Top\|\le 1$.
The function $f$ has no zeros in the unit disk, hence, 
the operators $\lambda-\Top$ are invertible for $\lambda\in\D$,
\begin{equation}
(\lambda-\Top)^{-1}(g+[f])=
\frac{g(z)-f(z)g(\lambda)/f(\lambda)}{\lambda-z}+[f],
\label{44f}
\end{equation}
the spectrum of $\Top$ lies on the unit circle,
\begin{gather}
\lim_{n\to+\infty}\|\Top^{-n}\|^{1/n}\le 1,\notag\\
\intertext{and}
(\lambda-\Top)^{-1}=-\sum_{n>0}\lambda^{n-1}\Top^{-n},
\qquad \lambda\in\D.
\label{44b}
\end{gather}

Next, we produce a (uniquely defined) continuous linear
operator $\Lop:\hil\longrightarrow
\ber/[f]$ coinciding with the canonical projection
$x\longmapsto x+[f]$
on $\ber$ and such that $\Top\circ\Lop=\Lop\circ \Sop$. Then,
$$\calM\overset{\text{def}}= [f]\subset \ker\Lop.$$
Moreover, $\ker\Lop$ is a nontrivial bilaterally invariant subspace
of $\hil$. Indeed,
$$
\Lop g=0\iff \Top\Lop g=0 \iff \Lop \Sop g=0. 
$$
By the definition of $\Lop$, we have $\ker\Lop\cap\ber=[f]$.
It follows that $f$ is non-cyclic in $\hil$ with respect to the
bilateral shift. Since $z^n-\Top^n(1+[f])\subset \calM_\T$, $n\in\Z$,
for every $x\in\hil$, we obtain $x-\Lop x\subset \calM_\T$, hence $x\in
\ber+\calM_\T$. Thus, $\ker\Lop=\calM_\T$ and $\hil=\ber+\calM_\T$.

The operator $\Lop$ is already defined on $\ber$. Furthermore, we
have $\Lop(z^{-n})=\Top^{-n}(1+[f])$, $n\in\Z$, and $\Lop$ extends by
linearity to the linear span of $\{z^{-n}\}_{n=1}^{+\infty}$.
We are to verify that $\Lop$ extends continuously to $\hil\ominus\ber$. For
every polynomial $q(z)=\sum_{n=1}^{+\infty}a_nz^{-n}$ in the variable
$z^{-1}$, we have
$$\Lop(q)=\sum_{n=1}^{+\infty}a_n\Top^{-n}(1+[f]),$$
the norm of which we estimate as follows:
\begin{multline*}
\|\Lop(q)\|_{\ber/[f]}\le
\sum_{n=1}^{+\infty}|a_n|\big\|\Top^{-n}(1+[f])\big\|_{\ber/[f]}\\
\le\biggl\{\sum_{n=1}^{+\infty}\frac{|a_n|^2}{\w(n-1)}\biggr\}^{1/2}
\biggl\{\sum_{n=1}^{+\infty}\big\|\Top^{-n}(1+[f])\big\|^2_{\ber/[f]}
\w(n-1)\biggr\}^{1/2}\\
=\|q\|_{\hil}
\biggl\{\sum_{n=1}^{+\infty}\big\|\Top^{-n}(1+[f])\big\|^2_{\ber/[f]}
\w(n-1)\biggr\}^{1/2}.
\end{multline*}
As a consequence, we have
$$
\big\|\Lop{\bigm|}_{\hil\ominus\ber}\big\|^2\le
\sum_{n=1}^{+\infty}\big\|\Top^{-n}(1+[f])\big\|^2_{\ber/[f]}\w(n-1).
$$
To estimate the right hand side of this inequality, we apply the 
identity
\begin{multline*}
A\overset{\text{def}}=\int_\D
\big\|(\lambda-\Top)^{-1}(1+[f])\big\|^2_{\ber/[f]}
\bas(\lambda)\mea(\lambda)\\=
\int_\D\bigg\|\sum_{n=1}^{+\infty}
\lambda^{n-1}\Top^{-n}(1+[f])\bigg\|^2_{\ber/[f]}
\bas(\lambda)\mea(\lambda)\\=
\sum_{n=1}^{+\infty}\big\|\Top^{-n}(1+[f])\big\|^2_{\ber/[f]}
\int_\D|\lambda|^{2n-2}\bas(\lambda)\mea(\lambda)\\=
\sum_{n=1}^{+\infty}\big\|\Top^{-n}(1+[f])\big\|^2_{\ber/[f]}\w(n-1),
\end{multline*}
which follows from \f{10a}, \f{44b}, and the fact that for radial
$\bas$ and for integers $n\ne k$,
$$
\int_\D \lambda^n\bar{\lambda}
{\vphantom\lambda}^k\bas(\lambda)\mea(\lambda)=0.
$$
Thus, to prove the theorem, we need only to verify that
$A<+\infty$. To estimate the norm of $(\lambda-\Top)^{-1}$ does not
seem very promising; fortunately, we do not need this. By
\f{44f},
\begin{gather*}
(\lambda-\Top)^{-1}(1+[f])=\frac{1-f(z)/f(\lambda)}{\lambda-z}+[f],\\
\intertext{hence we have}
\big\|(\lambda-\Top)^{-1}(1+[f])\big\|_{\ber/[f]}\le
\biggl\|\frac{1-f(z)/f(\lambda)}{\lambda-z}\biggr\|_{\ber},
\end{gather*}
and finally, we get 
\begin{gather*}
A\le \int_\D
\biggl\|\frac{1-f(z)/f(\lambda)}{\lambda-z}\biggr\|^2_{\ber}
\bas(\lambda)\mea(\lambda)\qquad\qquad\qquad\\
\qquad=\int_\D\int_\D \frac1{|f(\lambda)|^2}
\frac{|f(\lambda)-f(z)|^2}{|\lambda-z|^2}
\bas(\lambda)\bas(z)\mea(\lambda)\mea(z)<+\infty,
\end{gather*}
due to inequality \f{44a}. 
\qed
\bigskip
\bigskip

\section{\label{phrlind-est}A Phragm\'en--Lindel\"of type estimate for
functions in the disk}

Here, we will prove that if an analytic function $f$ is bounded in the unit
disk, and satisfies a kind of $\ell^p$ bound on a sequence of points
tending to the unit circle along a collection of circles, then we can bound
the analytic function in modulus by a given radial function, and the bound
is independent of the $H^\infty$-norm of $f$.

Later on, we will need the {\sl pseudo-hyperbolic metric} for the unit
disk, as given by
$$\rho_\D(z,w)=\left|\frac{z-w}{1-z\bar w}\right|,\qquad (z,w)\in\D^2.$$

Fix a constant $0<\kappa<1$, and let $\{r_n\}_n$ be a sequence of
numbers in the interval $[\frac45,1)$ tending to $1$ rather quickly. For
every $n=1,2,3,\ldots$, let $N_n$ be the integer that satisfies
$$
N_n\le \frac{\kappa}{1-r_n}<N_n+1.
$$ 
For each integer $k$ with $0\le k< N_n$, set
$$w_{n,k}=r_n\,e^{2\pi ik/N_n},$$
and select by some process (which will be explained later on in Section
\ref{main-constr}) a point $z_{n,k}$ from each disk
$$
\calD_{n,k}=\big\{z\in\D:\,|z-w_{n,k}|<(1-r_n)^2\big\}.
$$
Finally, consider a discrete measure $\mu$ equal to the sum of
point masses of size $1/(n^2N_n)$ at the points $z_{n,k}$, $0\le k<
N_n$, $n=1,2,3,\ldots$. In this section, we show that for
$0<p<+\infty$, the set of so-called {\sl analytic bounded point
evaluations} for $P^p(\mu)$ -- the closure of the polynomials in
$L^p(\mu)$ -- coincides with $\D$, and, more generally, we
supply effective estimates for the constants $C_p(z)$ in the
inequality
$$
|f(z)|^p\le C_p(z)\sum_{n,k}\frac{|f(z_{n,k})|^p}{n^2 N_n},
\qquad f\in H^\infty.
$$
Since we use only information on a discrete set in $\D$, we
cannot apply the standard technique of subharmonic functions.
Instead, we make use of the Lagrange interpolation formula;
the same method was applied earlier in \cite{HV}.
Let $B_n$ be the finite Blaschke product
$$
B_n(z)=\prod_{k=0}^{N_n-1}\frac{z-z_{n,k}}{1-\bar{z}_{n,k}z};
$$
an explicit calculation reveals that for $j=0,1,2,\ldots,N_n-1$, we have
$$\big(1-|z_{n,j}|^2\big)\,|B'(z_{n,j})|=\prod_{k:k\neq j}
\left|\frac{z_{n,j}-z_{n,k}}{1-\bar{z}_{n,k}z_{n,j}}\right|=\prod_{k:k\neq j}
\rho_\D\big(z_{n,j},z_{n,k}\big).$$

\begin{lemma} For some positive constant $c(\kappa)$, depending only on
$\kappa$, we have 

\noindent{\rm (a)} \qquad $|B'_n(z_{n,k})|\ge c(\kappa)\,N_n$,
\qquad $0\le k<N_n$.

\noindent
If $0<r<1$, $\varepsilon>0$ are fixed, and $r_n$ is sufficiently
close to $1$, then

\noindent{\rm (b)} \qquad $\bigl|\log|B_n(z)|+\kappa\bigr|\le\varepsilon$,
\qquad $|z|\le r$.
\label{tl8}
\end{lemma}

\begin{proof} Let $A_n$ be the finite Blaschke product
$$
A_n(z)=\prod_k\frac{z-w_{n,k}}{1-\bar{w}_{n,k}z}
=\frac{z^{N_n}-r^{N_n}_n}{1-r^{N_n}_nz^{N_n}},
$$
which is quite analogous to $B_n$. In view of our assumptions on the numbers
$N_n$ and the finite sequence $\{w_{n,k}\}_k$, we have
\begin{gather}
c_1\,N_n\le\big|A'_n(w_{n,j})\big|,\qquad j=0,1,2,\ldots,N_n-1,
\label{41a}\\
\lim_{n\to+\infty}\big|A_n(z)\big|=
\lim_{n\to+\infty}r^{N_n}_n=e^{-\kappa},\qquad |z|<1,\label{41b}
\end{gather}
where $c_1=c_1(\kappa)$ is a positive constant. The function $\log|B_n/A_n|=
\log|B_n|-\log|A_n|$ equals the sum of the functions
$$s_{n,k}(z)=\log\rho_\D(z,z_{n,k})-\log\rho_\D(z,w_{n,k})$$
over $k=0,1,\ldots N_n-1$. On the complement in $\D$ of the pseudohyperbolic
circle of radius $\frac12$ centered at $w_{n,k}$, we have
$$|s_{n,k}(z)|\le c_2\,(1-r_n)$$
for some positive constant $c_2=c_2(\kappa)$, whereas on the circle $r\T$,
we have
$$|s_{n,k}(z)|=o(1-r_n),\quad\text{ as}\quad r_n\to1;$$
along the unit circle $\T$, on the other hand, $s_{n,k}=0$. Summing up
all the terms, we obtain, after an application of the maximum principle,
that 
\begin{equation}
\Big|\log{|B_n(z)|}-\log{|A_n(z)|}\Big|=
o(1),\quad\text{ as}\quad r_n\to1,\quad\text{ for}\,\,\, |z|\le r,
\label{41c-1}
\end{equation}
and that for each $j=0,1,2,\ldots,N_n-1$,
\begin{equation}
\Big|\log|B'_n(z_{n,j})|-\log|A'_n(w_{n,j})|\Big|<c_3,
\label{41c-2}
\end{equation}
with some positive constant $c_3=c_3(\kappa)$ independent of the radius
$r_n$. 
We get both assertions (a) and (b) from the estimates \f{41a}, \f{41b},
\f{41c-1}, and \f{41c-2}.
\end{proof}

Let card be the function that computes the number of points in a given set
(it stands for {\sl cardinality}). 

\begin{lemma} For $n=1,2,3,\ldots$, put $\calN_n=\{0,\ldots,N_n-1\}$,
and take a subset $\calN^*_n$ of $\calN_n$. Let
$$
\sigma_n=\frac{\card(\calN_n\setminus \calN^*_n)}{N_n},
$$ 
be the density of $\calN_n^*$ in $\calN_n$, and define 
\begin{equation}
B^*_n(z)=\prod_{k\in \calN^*_n}\frac{z-z_{n,k}}{1-\bar{z}_{n,k}z}.
\label{05}
\end{equation}
If $r$ and $\varepsilon$, with $0<r<1$ and $0<\varepsilon<+\infty$
are fixed, and $r_n$ is sufficiently close to $1$,
then 
$$
-\varepsilon\le
\log|B^*_n(z)|+\kappa\le \frac{3\kappa \sigma_n}{1-|z|}+\varepsilon,
\qquad |z|\le r.
$$
\label{tl9}
\end{lemma}

\begin{proof} The assertion for $\sigma_n=0$ follows from
Lemma~\ref{tl8}(b). It remains to note that for $|z|\in[0,r]$ and 
for $r_n$ sufficiently close to $1$,
\begin{gather*}  
0\le\log|B^*_n(z)|-\log|B_n(z)|=
-\sum_{k\in \calN_n\setminus \calN^*_n}\log\rho_\D(z,z_{n,k})\\
\le\log\frac1{\rho_\D(|z|,r_n)}\,\card(\calN_n\setminus \calN^*_n)\le
\frac{3\kappa\sigma_n}{1-|z|}.
\end{gather*}
The proof is complete.
\end{proof}

\begin{prop} Suppose that the radii $r_n$ tend to $1$ sufficiently rapidly.
If $f$ is a function that is bounded and analytic in the unit disk, and
if for some $0<p<+\infty$,
\begin{equation}
\sum_{k=0}^{N_n-1}|f(z_{n,k})|^p\le n^2N_n,\qquad n=1,2,3,\ldots,
\label{03}
\end{equation}
then, for some constant $c=c(\kappa,p)$ independent of $f$,
$$
|f(z)|\le c\,\exp\left(\frac 1{1-|z|}\right),\qquad z\in\D.
$$
\label{tl10}
\end{prop}

\begin{proof} For each $n=1,2,3,\ldots$, we make the choice
$$
\calN^*_n=\Big\{k\in \calN_n:\,|f(z_{n,k})|\le n^{4/p}\Big\}.
$$
Let the associated finite Blaschke product $B_n^*$ be given by \f{05}. 
If $r_n$ tends to $1$ sufficiently rapidly, then by Lemma~\ref{tl9} and
the limit relationship \f{41b}, we have
\begin{equation}
-\frac 1{n^2}\le \log|B^*_n(z)|-N_n\log r_n\le 
\frac{4\kappa \sigma_n}{1-|z|}+\frac 1{n^2}, \qquad |z|\le r_{n-1}.
\label{081}
\end{equation}
Using \f{03} and a weak-type estimate, we obtain
\begin{equation}
\sigma_n=\frac{\card(\calN_n\setminus \calN^*_n)}{N_n}\le\frac 1{n^2}.
\label{082}
\end{equation}
Consider the set
$$\Gamma^*=\Big\{z_{n,k}:k\in \calN^*_n,\,\,\,n=1,2,3,\ldots\Big\},$$
which is a discrete lattice-like subset of $\D$. 
In view of \f{081} and \f{082}, the infinite product
$$
B^*(z)=\prod_{n=1}^{+\infty} \Big(r^{-N_n}_n\,B^*_n(z)\Big),
$$
converges uniformly on compact subsets of $\D$. Note that by \f{41b},
$r^{-N_n}_n\to e^{\kappa}$ as $n\to+\infty$. 
If we use the obvious estimate that $|B^*_n|<1$ on $\D$ for small values of
$n$, and the estimates \f{081} and \f{082} for large values of $n$, we easily
establish that for some positive constant $c_1=c_1(\kappa)$ that depends
only on $\kappa$, we have
\begin{equation}
|B^*(z)|\le c_1\,\exp\bigg(\frac1{5\,(1-|z|)}\bigg),\qquad z\in\D.
\label{06}
\end{equation}
Moreover, Lemma~\ref{tl8} -- modified to apply to $B^*_n$ instead of $B_n$,
using that $|B_n|\le|B^*_n|$ -- shows that if the radii $r_n$ tend to $1$
sufficiently rapidly as $n\to+\infty$, then
\begin{gather}
|B^*(z)|\ge c_2\,e^{\kappa n},\quad\text{for}\quad
|z|=1-2(1-r_n),\label{07}\\
\big|(B^*)'(z_{n,k})\big|\ge c_3\,N_n\,e^{\kappa n},\quad\text{for}\quad
k\in\calN_n^*,\,\,n=1,2,3,\ldots,
\label{08}
\end{gather}
where $c_2=c_2(\kappa)$ and $c_3=c_3(\kappa)$ are two positive constants.
By the Cauchy residue theorem, if $0<r<1$ and 
$r\T\cap\Gamma^*=\emptyset$, then
\begin{equation*}
-\frac{f(z)}{B^*(z)}+\sum_{z_{n,k}\in\Gamma^*\cap r\D}
\frac{f(z_{n,k})}{(B^*)'(z_{n,k})(z-z_{n,k})}
=\frac1{2\pi i}\int_{r\T}\frac{f(\zeta)}{B^*(\zeta)(z-\zeta)}
d\zeta,\qquad |z|<r.
\end{equation*}
We use the estimate \f{07} on circles $r\T$ and let $r\to1$, and realize that
the residue integral on the right hand side then tends to $0$, as it is
given that $f$ is bounded in $\D$. It follows that
$$\frac{f(z)}{B^*(z)}=\sum_{z_{n,k}\in\Gamma^*}
\frac{f(z_{n,k})}{(B^*)'(z_{n,k})(z-z_{n,k})},\qquad z\in\D,$$
where it is conceivable that the convergence is conditional. Taking absolute
values, we arrive at
$$
\Bigl|\frac{f(z)}{B^*(z)}\Bigr|\le\sum_{z_{n,k}\in\Gamma^*}
\frac{|f(z_{n,k})|}{\big|(B^*)'(z_{n,k})(z-z_{n,k})\big|},\qquad z\in\D.
$$
Hence, by \f{08}, for some positive constant $c_4=c_4(\kappa,p)$ independent
of $z$ and $f$, we have
$$
\Bigl|\frac{f(z)}{B^*(z)}\Bigr|\le \frac {c_4}{\dist(z,\Gamma^*)},
$$
where $\dist$ is the usual Euclidean distance function. If we make a
rather crude estimate of the above right hand side, and multiply by $|B^*|$,
then, by \f{06}, for some positive constant $c_5=c_5(\kappa,p)$,
we get
$$
|f(z)|\le c_5\,\exp\left(\frac1{4(1-|z|)}\right),
$$
provided that
$$|z|\in[0,1)\setminus\bigcup_n\left[1-2(1-r_n),1-\frac12\,(1-r_n)\right].$$
An application of the maximum principle in each ``complementary'' annulus
$$|z|\in\left[1-2(1-r_n),1-\frac12\,(1-r_n)\right]$$
completes the proof.
\end{proof}
\bigskip
\bigskip

\section{\label{aux-convex}The construction of harmonic functions: tools
from Convex Analysis}

In this section we start with the weight $\bas$, pass to an
associated function $\blm$ on the positive half-line and then
regularize it.

Suppose that a decreasing radial weight $\bas$ satisfies \rf{10},
and put
$$
\bdm(s)=\log\frac1{\bas(1-s)},\qquad 0<s<1.
$$
We want to study the behavior of $\bas(t)$ for $t$ near $1$, and
hence that of $\bdm(s)$ for $s$ near $0$. In order to study that
behavior in detail, we make an exponential change of coordinates,
and set
$$
\blm(x)=\log\bdm\big(e^{-x}\big)=\log\log\frac1{\bas(1-e^{-x})},
\qquad 0\le x<+\infty;
$$
the parameter
$$
x=\log\frac1{1-t}
$$
equals approximately the hyperbolic distance in $\D$ from $0$ to
$t$. The function $\blm(x)$ is increasing in $x$, and it grows
at least exponentially:
\begin{equation}
\lim_{x\to+\infty}e^{-\e x}\,\blm(x)=+\infty,
\label{90}
\end{equation}
where $\e>0$ is the same parameter as in \f{10}. We need the
following lemma from Convex Analysis.

\begin{lemma} Let $\blm:[0,+\infty)\longrightarrow(0,+\infty)$ be a 
$C^1$-smooth 
increasing function satisfying \rf{90}
for some $\e$, $0<\e\le1$. Then there exists another
$C^1$-smooth in\-creas\-ing convex function $\rlm:[0,+\infty)
\longrightarrow(0,+\infty)$, 
such that
\medskip

\noindent{\rm(a)} $\rlm$ is a minorant of $\blm$: $\rlm(x)\le\blm(x)$
holds on $[0,+\infty)$,
\smallskip

\noindent{\rm(b)} $e^{\e x}\le\rlm(x)$ holds on some 
interval $[A,+\infty)$, with $0\le A<+\infty$,
\smallskip

\noindent{\rm(c)} $\rlm'(x)\le\rlm(x)^{3/2}$ holds on $[0,+\infty)$,
\smallskip

and there is a sequence of numbers $\{x_n\}_n$ tending to $+\infty$, 
such that
\smallskip

\noindent{\rm(d)} $\rlm(x_n)=\blm(x_n)$, for all $n$,
\smallskip

\noindent{\rm(e)} $\e\,\rlm(x_n)\le\rlm'(x_n)$, for all $n$, and
\smallskip

\noindent{\rm(f)} $\rlm(x_n)+(x-x_n)\rlm'(x_n)+
\frac14\,\e\,(x-x_n)^2\,\rlm'(x_n)\le\rlm(x)$ for $x\in [0,+\infty)$.
\label{tl1}
\end{lemma}

\begin{proof} Put $Q=\sqrt{\blm}$, and write $q=\sqrt{\lambda}$, 
where $\lambda$ is the function we seek. The conditions (a)--(e)
then correspond to
\medskip

\noindent{\rm(a${}'$)} $q$ is a minorant of $Q$ on $[0,+\infty)$,

\noindent{\rm(b${}'$)} $e^{\e x/2}\le q(x)$ holds on $[A,+\infty)$,

\noindent{\rm(c${}'$)} $q'(x)\le\frac12\,q(x)^2$ on $[0,+\infty)$,

\noindent{\rm(d${}'$)} $q(x_n)=Q(x_n)$ for all $n$,

\noindent{\rm(e${}'$)} $\frac12\,\e\,q(x_n)\le q'(x_n)$ for all $n$.
\medskip

We shall find a convex $q$ with the above properties, which
means that $\lambda=q^2$ is convex as well. Note that by the
assumption on $\Lambda$, 
\begin{equation}
\lim_{x\to+\infty}e^{-\e x/2}\,Q(x)=+\infty,
\label{24br}
\end{equation}
so that if we forget about property (c${}'$), we can just take $q$ to
be equal to the greatest convex minorant $q_0$ of $Q$; it is not hard
to check
that all the properties (a${}'$), (b${}'$), (d${}'$), and (e${}'$) 
are fulfilled for a sequence $\{x_n\}_n$.
To get also (c${}'$), we apply an iterative procedure.

As a consequence of \f{24br}, we get
\begin{equation}
\lim_{x\to+\infty}e^{-\e x/2}q_0(x)=+\infty.
\label{122}
\end{equation}
This is so because the function $q_0$ must touch $Q$ along an
unbounded closed set, which we denote by $\calE$. Note that $q_0$ is
affine on each open interval in the complement of $\calE$. 
Changing $q_0$ a little
on a small interval with the origin as the left end point, if
necessary, we can
guarantee that $q'_0(0)\le\frac13\,q_0(0)^2$. 
Set $a_0=0$. Our iterative procedure runs as follows. For every
$k=0,1,2,3,\ldots$, we start with a convex increasing minorant
$q_k$ of $Q$ such that $q'_k(a_k)\le\frac13\,q_k(a_k)^2$, 
$q'_k(x)\le\frac12\,q_k(x)^2$ for $x\le a_k$, and
$q_k(x)=q_0(x)$ for $x\ge a_k$.
If
$$
q'_k(x)\le\frac12\,q_k(x)^2
$$
on the whole interval $[a_k,+\infty)$, then we are done, because
we pick $q=q_k$. If it is not so, there exists of course a
point $x\in(a_k,+\infty)$ with
$$
q'_k(x)>\frac12\,q_k(x)^2.
$$
Let $b_k\ge a_k$ be the infimum of all such points $x$. 
Let $(c_k,d_k)$ be the maximal interval in $(a_k,+\infty)$ 
such that $b_k\in(c_k,d_k)$ and
$$
q'_k(x)>\frac13\,q_k(x)^2,\qquad x\in(c_k,d_k).
$$
Notice that $q'_k(c_k)=\frac13\,q_k(c_k)^2$. We claim that
$c_k$ belongs to the set $\calE$ we defined earlier. In fact, every
small interval to the right of $c_k$ contains infinitely many
points of $\calE$. The reason is that away from $\calE$, $q_k$ is
affine, and since $q_k$ is increasing, we get
$$
\frac{q_k(c_k)^2}{3}\le\frac{q_k(x)^2}{3}<q'_k(x)=q'_k(c_k)
=\frac{q_k(c_k)^2}3
$$
if for some $x$ with $c_k<x<d_k$, we have $(c_k,x)\cap \calE=\emptyset$; 
this is impossible. We shall now alter the function $q_k$ to the right
of the point $c_k$. 
Let $f_k$ solve the initial value problem
$$
f'_k(x)=\frac13\,f_k(x)^2, \qquad f_k(c_k)=q_k(c_k);
$$
we see that $f_k$ explodes in finite time. As a matter of fact,
an explicit calculation reveals that
$$
f_k(x)=\frac3{c_k+\frac3{q_k(c_k)}-x},\qquad x\in I_k,
$$
where $I_k=[c_k,c_k+\hbox{$\frac3{q_k(c_k)}$})$, and the
explosion point is the right end point of the indicated
interval. Now, take as $q_{k+1}$ the greatest convex minorant
of the function that equals $q_k$ on $[0,+\infty)\setminus I_k$
and equals $\min\{f_k,q_k\}$ on $I_k$. The function $f_k$ is
convex and increasing on $I_k$, and to the right of $c_k$, it
initially drops below the convex function $q_k$, but then after
a while, it grows above it again, finally to explode at the
right end point of $I_k$. It follows that $q_{k+1}$ is
increasing, $q_{k+1}=q_k$ on the interval $[0,c_k]$,
$q_{k+1}=f_k$ on some interval $[c_k,e_k]\subset I_k$, $q_{k+1}$
is affine on some interval $[e_k,a_{k+1}]$, at the right end
point of which $q_{k+1}$ touches $q_k$, and $q_{k+1}=q_k$ on the
interval $[a_{k+1},+\infty)$. Since $q_{k+1}$ is increasing,
we have
$$
\frac{q'_{k+1}(x)}{q_{k+1}(x)^2}=\frac{q'_{k+1}(e_k)}
{q_{k+1}(x)^2}\le
\frac{q'_{k+1}(e_k)}{q_{k+1}(e_k)^2}=
\frac{f'_k(e_k)}{f_k(e_k)^2}=\frac 13,
\qquad x\in [e_k,a_{k+1}].
$$
Hence, 
\begin{equation}
a_{k+1}\ge d_k.
\label{21br}
\end{equation}
We get that $q_{k+1}$ is a convex in\-creas\-ing minorant of $Q$,
with $$q'_{k+1}(a_{k+1})\le\frac13\, q_{k+1}(a_{k+1})^2,$$
and 
$$q'_{k+1}(x)\le\frac12\,q_{k+1}(x)^2\,\,\hbox{ for }\,x\le a_{k+1},
\quad q_{k+1}(x)=q_0(x)\,\,\hbox{ for }\,x\ge a_{k+1}.$$

If our iterative procedure does not stop on a finite step (that
is, if $q'_k(x)>\frac12\,q_k(x)^2$ on an unbounded subset of the
interval $[0,+\infty)$), then we put $q=\lim_{k\to\infty}q_k$.
Next we verify the properties (a${}'$)--(e${}'$) for $q$.

Notice that $\lim a_k=\infty$. Indeed, by the definition of
$b_k$, $(c_k,d_k)$, and by the property \f{21br} of $a_k$, 
we have $a_k\le c_k<b_k<d_k\le a_{k+1}$. 
If $a_k$ were to tend to a finite number $a_\infty$ as
$k\to+\infty$, then $c_k\to a_\infty$, $b_k\to a_\infty$, as
$k\to+\infty$ as well. By the definitions of these points,
$$
q'_0(c_k)=\frac13\,q_0(c_k)^2,\qquad q'_0(b_k)
=\frac12\,q_0(b_k)^2.
$$
In the limit, we obtain a contradiction, which does it.

It follows from our construction that $q=q_k$ on $[0,a_k]$,
$k=1,2,3,\ldots$, and that $q$ is a convex minorant of $Q$.
Hence, (a${}'$) and (c${}'$) follow.
The graph of function $q$ touches that of $q_0$ except on the intervals
$(c_k,a_{k+1})$, so that if $q$ does not grow too slowly on
those intervals, we have
$$
\lim_{x\to+\infty}e^{-\e x/2}q(x)=+\infty.
$$
We now show that the above holds in general. We know that
the function
$e^{-\e x/2}q_{k+1}(x)$ is sufficiently big for
$x\in\{c_k,a_{k+1}\}$. 
Also, the function
$$
e^{-\e x/2}f_k(x)=e^{-\e x/2}q_{k+1}(x),\qquad x\in[c_k,e_k],
$$
increases on $[c_k,e_k]$ (at least for big $k$). 
Since $q_k$ is affine on $[e_k,a_{k+1}]$, we obtain
$q(x+e_k)=q_k(e_k)+x\,q'_k(e_k)$, $x\in [0,a_{k+1}-e_k]$.
Furthermore, a simple calculation yields that the function
$$
t\mapsto e^{-\e(t+e_k)/2}\big[q_k(e_k)+t\,q'_k(e_k)\big]
$$
either is monotonic or has just one local maximum on $[0,+\infty)$.
As a consequence,
$$
e^{-\e x/2}q_k(x)\ge\min\big\{e^{-\e x/2}q_k(x):
\,x\in\{c_k,a_{k+1}\}\big\},
\qquad e_k\le x\le a_{k+1},
$$
which finishes the proof of (b${}'$).

We turn to verifying that there exists a sequence $\{x_n\}_n$ of
points in $\calE$ tending to $+\infty$ for which
\begin{equation}
q(x_n)=q_0(x_n),\qquad \frac12\,\e\,q(x_n)<q'(x_n).
\label{123}
\end{equation}
The points $c_k$ are in $\calE$, and they may work as $x_n$,
provided that there are infinitely many of them. If there is
only a finite supply of $c_k$, then $q(x)=q_0(x)$ for all
sufficiently big $x$. If we then also cannot find arbitrarily
big $x_n\in \calE$ satisfying \f{123}, then the reason is that
$q'_0(x)\le\frac12\e q_0(x)$ for all $x$ in $\calE$ that are
sufficiently large, and since $q_0$ is affine outside $\calE$, we
get a contradiction with \f{122}. Thus, (d${}'$) and (e${}'$)
follow.

As we made clear before, we pick $\rlm(t)=q(t)^2$. We need only
verify property (f). Let $x_n$ be as in \f{123}, and take an
arbitrary $x\in[0,+\infty)$. We then have
$$
q(x_n)+(x-x_n)q'(x_n)\le q(x)
$$
and hence
\begin{align*}
\rlm(x)&\ge q(x_n)^2+2(x-x_n)q(x_n)q'(x_n)+(x-x_n)^2q'(x_n)^2\\
&=\rlm(x_n)+(x-x_n)\rlm'(x_n)+(x-x_n)^2\rlm'(x_n)
\frac{q'(x_n)}{2q(x_n)}\\
&\ge\rlm(x_n)+(x-x_n)\rlm'(x_n)+\frac14\,\e\,(x-x_n)^2\rlm'(x_n). 
\end{align*}
The proof is complete.
\end{proof}
\bigskip
\bigskip

\section{\label{aux-harm-constr}The construction of harmonic functions:
building blocks}

We start with the weight $\bas$, pass on to the associated functions $\bdm$
and $\blm$ on the interval $(0,1]$ and the positive half-line,
respectively, and then apply Lemma~\ref{tl1} to $\blm$ in order to obtain the
minorant $\rlm$ to $\Lambda$ along with the sequence $\{x_n\}_n$. We recall
the relationships
$$\bdm(s)=\log\frac1{\bas(1-s)}\quad\text{ and }\quad
\blm(x)=\log\bdm(e^{-x}),$$
and put 
\begin{equation} 
\rdm(s)=\exp\Bigl[\rlm\Bigl(\log\frac 1s\Bigr)\Bigr],\qquad 0<s<1.
\label{01}
\end{equation} 
Then $\exp(s^{-\e})\le\rdm(s)\le\bdm(s)$, for small $s>0$, and we
also have
$\rdm(e^{-x_n})=\bdm(e^{-x_n})$; as usual, $\e$ is the positive
quantity which appears in \f{10}. 

For the rest of this section, we fix some sufficiently big
$x_n$ from the above-mentioned sequence, and put
$$\delta_n=e^{-x_n},\quad r_n=1-\delta_n, \quad\text{ and}\quad
\g_n=e^{-\rlm(x_n)/10}.$$
Furthermore, let
\begin{equation}
\left.
\begin{gathered}
H_n(z)=e^{\rlm(x_n)-x_n\rlm'(x_n)}\,(1-z)^{-\rlm'(x_n)}\\
=e^{\rlm(x_n)-x_n\rlm'(x_n)}\,
\exp\bigg[\rlm'(x_n)\log\frac1{1-z}\bigg],\qquad z\in\D,\\
h_n=\re H_n.
\end{gathered}
\right\}
\label{0101}
\end{equation}
The function $H_n$ is holomorphic on $\D$, and consequently, $h_n$
is real-valued and harmonic. There are sectors about the point $1$
of angular opening $\pi/\rlm'(x_n)$ where $h_n$ is large, positive
or negative, in an alternating fashion. We shall use the functions
$h_n$ as ``building blocks'', like in \cite{BH}. The intention is to
construct a function that is alternatingly big and small near a
prescribed sequence of points on the unit circle.

For the construction, we need sufficiently effective estimates on the
growth and decay of $h_n$.

In the following lemma, we compare the size of $h_n(z)$ with that of
the majorants 
$$M(z)=\rdm(1-|z|),\quad M_n(z)=(1+\g_n)^{-1}
\bigl[\rdm(1-|z|)-2[\log\rdm(1-|z|)]^2\bigr];$$
clearly, $M_n\le M$. It is easy to see that with our choice of $h_n$ in 
formula (\ref{0101}), and the clever choice of the point $x_n$ of Lemma
\ref{tl1} we have $|h_n|\le M$ throughout $\D$, and 
$h_n(r_n)=M(r_n)$.
We prove that $h_n\le M_n$ holds outside a small disk $\calD_n$ of radius 
$(1-r_n)^2$ centered at $r_n$, and that $-h_n\le M_n$ in a neighborhood 
of the unit circle that does not depend on the point $x_n$.

\begin{lemma} In the above setting, we have, for $z\in\D$,
\smallskip

\noindent{\rm (a)}\quad $|h_n(z)|\le\rdm(1-|z|)$,
\smallskip

\noindent{\rm (b)}\quad
$(1+\g_n)\,h_n(z)\le\rdm(1-|z|)-2[\log\rdm(1-|z|)]^2$ when 
$|z-r_n|>\delta^2_n$, $\,$and 
\smallskip

\noindent{\rm (c)}\quad 
$(1+\g_n)\,h_n(z)>-\rdm(1-|z|)+2[\log\rdm(1-|z|)]^2\,\,\, 
\text{when}\,\,\,c(\rdm)<|z|<1$, where $c(\rdm)$, $0<c(\rdm)<1$, only
depends on the function $\rdm$.
\label{tl4}
\end{lemma}

The next step is to verify the following additional properties of the
functions $h_n$: 
\medskip

\noindent$\bullet$ $\,\,|h_n|$ is small outside a disk of radius proportional 
to $1-r_n$ centered at $1$, 
\medskip

\noindent$\bullet$ $\,\,-h_n$ is sufficiently big at a point inside this 
disk, and
\medskip

\noindent$\bullet$ $\,\,(1+\g_n)\,h_n-\bdm(1-|z|)$ is sufficiently big 
(in the integral sense) inside the above-mentioned disk $\calD_n$.
\medskip

\begin{lemma} Under the conditions of Lemma~{\rm\ref{tl4}},
\smallskip

\noindent{\rm (a)}\quad 
$|h_n(z)|<\delta_n^3\quad \text{when}\quad z\in\D\quad\text{has}\quad
|z-1|\ge\delta_n\exp (-1+2/\e)$,
\smallskip

\noindent{\rm (b)}\quad 
$h_n(w)=-\rdm(\delta_n)\quad \text{holds for some $w\in\D$ 
with}\quad1-\delta_n<|w|<\frac12\,\delta_n$, and 
\smallskip

\noindent{\rm (c)}\quad we have the integral estimate
$$1\le\int_{\calD_n}\exp\big[(1+\g_n)h_n(z)-\bdm(1-|z|)\big]\mea(z).$$
\label{tl41}
\end{lemma}

For the proofs of these lemmas, we use the following simple estimates.

\begin{lemma} For positive $\alpha$, let $F_\alpha$ denote the function
$$
F_\alpha(z)=(1-z)^{-\alpha},\qquad z\in\D,
$$
where the power is defined by the principal branch of the logarithm. 
For $0<\varphi<\pi/6$, let $\calR_\varphi$ be the domain
$$
\calR_\varphi=\Big\{z\in\C:\,\tfrac12<|z|<1,\,\,
\varphi<|\arg(1-z)|\Big\},
$$
and $\seg_\alpha$ the union of two line segments
$$
\seg_\alpha=\Big\{z\in\D:\,|\arg(1-z)|=\frac{\pi}{\alpha}\Big\}.
$$
Then the real part of $F_\alpha$ has the following properties: 
\smallskip

\noindent{\rm(a)} 
$\re F_\alpha(z)=|1-z|^{-\alpha}\cos\big(\alpha\arg(1-z)\big),
\qquad z\in\D$,
\smallskip

\noindent{\rm(b)} 
$\re F_\alpha(z)\le|F_\alpha(z)|\le 
F_\alpha(|z|)=(1-|z|)^{-\alpha},\qquad
z\in\D$,
\smallskip

\noindent{\rm(c)} 
$\re F_\alpha(z)=-|1-z|^{-\alpha},\qquad z\in \seg_\alpha$,
\smallskip

and, if $\alpha$ is sufficiently large, 
\smallskip

\noindent{\rm(d)} 
$\re F_\alpha(z)\le |F_\alpha(z)|\le\exp
\Bigl(-\dfrac{\alpha\varphi^2}{3}\Bigr)\,
\big(1-|z|\big)^{-\alpha},
\qquad z\in\calR_\varphi$.
\label{tl2}
\end{lemma}

\begin{proof} Equality (a) follows from the definition of the
power function; equality (c) is an immediate consequence of (a).
The inequality (b) follows from the triangle inequality.

Finally, for large $\alpha$, a geometric consideration using the
inequality
$$
\cos\varphi\le 1-\frac{\varphi^2}3,\qquad 0<\varphi<\pi/6,
$$ 
yields the estimate
$$
1-|z|\le\Bigl(1-\frac{\varphi^2}3\Bigr)|1-z|,\qquad
z\in\calR_\varphi,
$$
and (d) follows.
\end{proof}

\begin{proofof}{Lemma \ref{tl4}} Fix $z\in\D$, 
put
$$t=\log \frac1{1-|z|},$$
and consider the affine function
$$
\aff_n(t)=\rlm(x_n)-x_n\rlm'(x_n)+t\rlm'(x_n).
$$ 
Then $\aff_n(x_n)=\lambda(x_n)$ and $\aff_n'(x_n)=\rlm'(x_n)$, so that
by the convexity of $\rlm$, $\aff_n(t)\le\rlm(t)$ holds everywhere. 
Using Lemma \ref{tl2}(b), we get
\begin{equation*}
\log|h_n(z)|\le\log|H_n(z)|
\le\rlm(x_n)-x_n\rlm'(x_n)+t\rlm'(x_n)\le
\rlm(t)=\log\rdm(1-|z|).
\end{equation*}
This proves part (a).
\smallskip

In the same fashion, we see that keeping the notation
$$t=\log\frac1{1-|z|},$$
we get
\begin{equation}
\big|\tilde h_n(z)\big|\le \exp\big[\aff_n(t)+\g_n\big]\le\exp[\rlm(t)+\g_n],
\label{74}
\end{equation}
where
$$\tilde h_n(z)=(1+\g_n)\,h_n(z),\qquad z\in\D.$$
Moreover, by Lemma \ref{tl1}(f),
\begin{multline}
\rdm(1-|z|)=\exp[\rlm(t)]
\ge\exp\left[\rlm(x_n)+(t-x_n)\rlm'(x_n)+
\frac14\,\e(t-x_n)^2\rlm'(x_n)\right]\\
=\exp\left[\aff_n(t)+\frac14\,\e(t-x_n)^2\rlm'(x_n)\right].
\label{71}
\end{multline}
It follows that
\begin{multline}
\rdm(1-|z|)-\big|\tilde h_n(z)\big|=\rdm(1-|z|)\,
\left(1-\frac{\big|\tilde h_n(z)\big|}{\rdm(1-|z|)}\right)\\
\ge\left(1-\exp\left[\g_n-\frac14\,\e(t-x_n)^2\rlm'(x_n)\right]\right)
\,\exp[\rlm(t)],
\label{72}
\end{multline}
and that
\begin{equation}
\rdm(1-|z|)-\big|\tilde h_n(z)\big|
\ge\left(\exp\left[\frac14\,\e(t-x_n)^2\rlm'(x_n)\right]-
e^{\g_n}\right)\,\exp[\aff_n(t)].
\label{75}
\end{equation}
\smallskip
If $z$ belongs to the domain under consideration in (b), then
$$\bigl| |z|-r_n\bigr|>\frac12\,\delta_n^2\quad\text{or}\quad
|\arg(1-z)|>\frac{\delta_n}2.$$
Furthermore, if $h(z)<0$, then $|\arg(1-z)|>\frac12\pi/\rlm'(x_n)$.
Write
$$\eta_n=\min\left\{\frac{\delta_n}2,\frac{\pi}{2\rlm'(x_n)}\right\}.$$
To verify (b) and (c), we need to estimate the expression
$$\rdm(1-|z|)-\big|\tilde h_n(z)\big|,$$
for points $z\in\D$ that satisfy at least one of the following two conditions:
\smallskip

\noindent (A) $\quad\bigl||z|-r_n\bigr|>\delta_n^2/2$ and

\noindent (B) $\quad|\arg(1-z)|>\eta_n$.
\smallskip

\noindent We first look at the case (A). Note that if
$$\bigl||z|-r_n\bigr|>\frac{\delta_n^2}2,$$
then
$$
|t-x_n|=\Bigl|\log\frac{1-r_n}{1-|z|}\Bigr|>\frac14\,\delta_n
=\frac14\,e^{-x_n}.
$$

\noindent (i) If $t-x_n>\frac14\,e^{-x_n}$, then, for big $x_n$,
$$
\frac{\e}4\,(t-x_n)^2\rlm'(x_n)>e^{-\rlm(x_n)/{3}}+\g_n,
$$
and by \f{72},
\begin{equation}
\rdm(1-|z|)-\big|\tilde h_n(z)\big|\ge\bigg(1-
\exp\left[-e^{-\rlm(x_n)/3}\right]\bigg)
\exp[\rlm(t)]\ge\exp\left(\frac{\rlm(t)}2\right). 
\label{129}
\end{equation}

\noindent (ii) If
$$\frac14\,e^{-x_n}<x_n-t<\frac{\rlm(x_n)}{2\,\rlm'(x_n)},$$
then for big values of $x_n$,
$$
\frac14\,\e\,(t-x_n)^2\rlm'(x_n)>\exp\left(-\frac{\rlm(x_n)}5\right)+\g_n
$$
and
$$
\aff_n(t)>\frac12\,\rlm(x_n),
$$
so that by (\ref{75}),
\begin{equation}
\rdm(1-|z|)-\big|\tilde h_n(z)\big|\ge \exp\left(\frac{\rlm(x_n)}4\right).
\label{73}
\end{equation}
Next, suppose that $x_n-t\ge\frac12\,\rlm(x_n)/\rlm'(x_n)$. Then, by 
Lemma~\ref{tl1}(c),
\begin{equation}
\frac1{16}\,\e\,[\rlm(x_n)]^{1/2}\le
\frac{\e\,[\rlm(x_n)]^2}{16\,\rlm'(x_n)}
\le\frac14\,\e\,(t-x_n)^2\rlm'(x_n)
\label{78}
\end{equation}

\noindent (iii) Thus, if
$$\frac{\rlm(x_n)}{2\,\rlm'(x_n)}\le x_n-t<\frac{\rlm(x_n)}{\rlm'(x_n)},$$
then $\aff_n(t)>0$, and by (\ref{75}),
\begin{equation}
\rdm(1-|z|)-\big|\tilde h_n(z)\big|>
\exp\left(\frac1{32}\,\e\,\rlm(x_n)^{1/2}\right). 
\label{100}
\end{equation}

\noindent (iv) Finally, if
$$x_n-t\ge\frac{\rlm(x_n)}{\rlm'(x_n)},$$
then $\aff_n(t)\le 0$. By (\ref{74}),
$$\big|\tilde h_n(z)\big|\le e^{\g_n}\le e,$$
and
\begin{equation}
\rdm(1-|z|)-\big|\tilde h_n(z)\big|\ge\exp[\rlm(t)]-e.
\label{101}
\end{equation}
As a result of \f{129}, \f{73}, \f{100} and \f{101}, for
$$\big||z|-r_n\big|>\frac{\delta_n^2}2,$$
and for big $x_n$ and $t$, we obtain
\begin{equation}
\rdm(1-|z|)-\big|\tilde h_n(z)\big|>
\exp\left(\frac1{32}\,\e\,\rlm(t)^{1/2}\right)-e
\ge 2[\rlm(t)]^2=2[\log\rdm(1-|z|)]^2. 
\label{102}
\end{equation}
\smallskip

\noindent (B) In this case,
$$|\arg(1-z)|>\eta_n=\min\left\{\frac{\delta_n}2,\,
\frac{\pi}{2\,\rlm'(x_n)}\right\},$$ 
and by Lemma~\ref{tl2}(d), we have
$$
\big|\tilde h_n(z)\big|\le\exp\left(\aff_n(t)+\g_n-\frac13\,
\rlm'(x_n)\,\eta_n^2\right).
$$
By Lemma~\ref{tl1}(c),
$$
\g_n-\frac{\rlm'(x_n)\,\eta_n^2}3<-\frac{\rlm'(x_n)\,\eta_n^2}4.
$$
Hence,
\begin{equation}
\big|\tilde h_n(z)\big|\le
\exp\left(\aff_n(t)-\frac14\,\rlm'(x_n)\eta_n^2\right).
\label{83}
\end{equation}
Using \f{71}, we get, analogously to \f{72} and \f{75},
that
\begin{multline}
\rdm(1-|z|)-\big|\tilde h_n(z)\big|\\
\ge\left(1-\exp\left[-\frac14\,\rlm'(x_n)\eta_n^2-\frac14\,
\e\,(t-x_n)^2\rlm'(x_n)\right]\right)\,\exp[\rlm(t)],
\label{76}
\end{multline}
and
\begin{multline}
\rdm(1-|z|)-\big|\tilde h_n(z)\big|\\
\ge\left(\exp\left[\frac14\,\e\,(t-x_n)^2\rlm'(x_n)\right]-
\exp\left[-\frac14\,\rlm'(x_n)\eta_n^2\right]\right)\,\exp[\aff_n(t)].
\label{77}
\end{multline}

\noindent (i) If $t>x_n$, then by \f{76} and by Lemma~\ref{tl1}, parts
(c) and (e), we have, for big $x_n$,
\begin{equation}
\rdm(1-|z|)-\big|\tilde h_n(z)\big|>\exp\left(\frac{\rlm(t)}2\right). 
\label{103}
\end{equation}

\noindent (ii) If
$$0\le x_n-t<\frac{\rlm(x_n)}{2\,\rlm'(x_n)},$$
then $\aff_n(t)>\rlm(x_n)/2$, and by \f{77} and by Lemma~\ref{tl1}, parts
(c) and (e), we have, for big $x_n$,
\begin{equation}
\rdm(1-|z|)-\big|\tilde h_n(z)\big|>\exp\left(\frac13\,\rlm(x_n)\right).
\label{103a}
\end{equation}
If
$$x_n-t\ge\frac{\rlm(x_n)}{2\,\rlm'(x_n)},$$
then we argue as in case (A) using the estimate \f{78}.

\noindent (iii) If
$$\frac{\rlm(x_n)}{2\,\rlm'(x_n)}\le x_n-t
<\frac{\rlm(x_n)}{\rlm'(x_n)},$$
then $\aff_n(t)>0$, and by (\ref{77}),
\begin{equation}
\rdm(1-|z|)-\big|\tilde h_n(z)\big|>\exp\left(\frac1{32}\,
\e\,[\rlm(x_n)]^{1/2}\right). 
\label{79}
\end{equation}

\noindent (iv) Finally, if $x_n-t\ge\rlm(x_n)/\rlm'(x_n)$, then
$\aff_n(t)\le 0$.
By (\ref{83}), $\big|\tilde h_n(z)\big|\le 1$, and 
\begin{equation}
\rdm(1-|z|)-\big|\tilde h_n(z)\big|>\exp[\rlm(t)]-1.
\label{80}
\end{equation}

As a result of \f{103}--\f{80}, for $|\arg(1-z)|>\eta_n$, 
and for big $x_n$ and $t$, we get
\begin{equation}
\rdm(1-|z|)-\big|\tilde h_n(z)\big|>\exp
\left[\frac1{32}\,\e\,[\rlm(t)]^{1/2}\right]-1
\ge 2\,[\rlm(t)]^2=2[\log\rdm(1-|z|)]^2. 
\label{104}
\end{equation}
\smallskip

The estimates \f{102} and \f{104} imply both (b) and $(c)$ for big
values of $x_n$.
\end{proofof}

\begin{proofof}{Lemma \ref{tl41}} To verify (a), we
note that if
$$|z-1|\ge\delta_n\,e^{2/\e},\quad\text{ then}\quad
\log\frac1{|1-z|}<x_n-\frac2\e,$$
and by Lemma \ref{tl1}, parts (b) and (e), together with Lemma
\ref{tl2}(b), we have, for big $x_n$,
\begin{multline*}
\log|h_n(z)|\le \rlm(x_n)-x_n\rlm'(x_n)+x_n\rlm'(x_n)-
\left(\frac{2}{\e}-1\right)\,
\rlm'(x_n)=\\
\rlm(x_n)-\left(\frac{2}{\e}-1\right)\,\rlm'(x_n)<-3\,x_n,
\end{multline*}
as desired.
\smallskip

To prove (b), note that if
$$w=1-\delta_n\exp\left(i\,\frac{\pi}{\rlm'(x_n)}\right),$$
then we have $\frac12\,\delta_n<1-|w|<\delta_n$ and $w\in\seg_{\rlm'(x_n)}$,
so that in view of part (c) of Lemma~\ref{tl2}, we obtain
$$
h_n(w)=-\exp\big[\rlm(x_n)-x_n\rlm'(x_n)+x_n\rlm'(x_n)\big]=
-\rdm(\delta_n).
$$
\smallskip

Finally, to prove the estimate in part (c), we consider the region
$$
\calR_n=\Big\{r\,e^{i\theta}:\,r_n-\g_n^2<r<r_n,\,\,\,|\theta|<\g_n^2\Big\}.
$$
which is a subset of the domain of integration in (c).
Using that for $z\in\calR_n$,
\begin{gather*}
\delta_n\leq |1-z|\le \delta_n+2\gamma_n^2,\\
|\arg(1-z)|\le 2\,\frac{\gamma_n^2}{\delta_n},\\
\bdm(1-|z|)\le\bdm(\delta_n)=\rdm(\delta_n),
\end{gather*}
we derive from Lemma \ref{tl2}, part (a) and Lemma \ref{tl1}, parts (b)
and (c), together with the convexity of the function $\lambda$, that
for $z\in\calR_n$, 
\begin{multline*}
\log h_n(z)+\log(1+\g_n)=\rlm(x_n)-x_n\rlm'(x_n)\\
+\rlm'(x_n)\log\frac1{|1-z|}
+\log\Big(\cos\big[\rlm'(x_n)\arg(1-z)\big]\Big)+\log(1+\g_n)\\
\ge\rlm(x_n)-x_n\rlm'(x_n)+\rlm'(x_n)\log\frac1{\delta_n}+\frac{2\g_n}{3}
-\frac{4\g_n^2\rlm'(x_n)}{\delta_n}-\frac{4\g_n^4\rlm'(x_n)^2}{\delta_n^2}\\
\ge\rlm(x_n)-x_n\rlm'(x_n)+x_n\rlm'(x_n)+\frac{\g_n}{2}=
\log\bdm(\delta_n)+\frac{\g_n}{2}.
\end{multline*}
As a consequence, it follows that for large $x_n$,
\begin{gather*}
\int_{\calD_n}
\exp\big[(1+\g_n)\,h_n(z)-\bdm(1-|z|)\big]\mea(z)\\
\ge\int_{\calR_n}\exp\big[(1+\g_n)\,h_n(z)-\bdm(1-|z|)\big]\mea(z)\\
\ge m_\D(\calR_n)\exp\big[\bdm(\delta_n)e^{\g_n/2}-\bdm(\delta_n)\big]
\ge\g_n^4\exp\left[\frac12\,\gamma_n\bdm(\delta_n)\right]\\
=\exp\left[-\frac25\,\lambda(x_n)\right]\,\exp
\left[\frac12\,e^{-\lambda(x_n)/10}e^{\lambda(x_n)}\right]\ge 1.
\end{gather*}
This completes the proof of part (c), and hence that of the whole lemma.
\end{proofof}

\section{\label{aux-harm-est}The construction of harmonic functions:
estimates}

For the proof of Theorem 5.3, we need additional estimates on
the values of $h$ and $H'$ at pairs of nearby points in the unit disk.
We begin with a simple regularity lemma.

\begin{lemma} For positive $s$ close to $0$,
$$
\rdm\big(s-[\rdm(s)]^{-2}\big)<\rdm(s)+1.
$$
\label{tl5}
\end{lemma}

\begin{proof} We start with the inequality
$$
e^t\rlm'(t)<e^{\rlm(t)/2},
$$
which follows from Lemma~\ref{tl1}, parts (b) and (c). Using \f{01} and
passing to the variable $x=e^{-t}$, we get
$$
|\rdm'(x)|<\rdm(x)^{3/2}.
$$
Since $\rdm$ is monotonically decreasing, to prove the lemma, it suffices to
notice that if $t<s$ and $\rdm(t)=\rdm(s)+1$, then
$\rdm(s)\le\rdm(x)\le\rdm(s)+1$ for $t\le x\le s$, and hence
$$
1=\rdm(t)-\rdm(s)=-\int^s_t\rdm'(x)\,dx<\int^s_t[\rdm(x)]^{3/2}\,dx
\le(s-t)[\rdm(s)+1]^{3/2},
$$
so that $s-t>[\rdm(s)]^{-2}$. 
\end{proof}

In the following lemma, for points $z,w\in\D$ that are sufficiently close
to one another, we produce upper estimates for $h_n(w)$ and 
$2\log|H_n'(w)|+h_n(w)$ that depend on the size of $|z|$.

\begin{lemma} In the notation of Lemma~{\rm\ref{tl4}}, take
$w\in\D$ with $|w-z|<[\rdm(1-|z|)]^{-2}$ for some $z\in\D$
with $|z|$ sufficiently close to $1$. Then
\begin{equation}
\big|\rdm(1-|z|)-\rdm(1-|w|)\big|<1.
\label{41}
\end{equation}

\noindent{\rm (a)} If $\rlm(t)\le\frac25\,\rlm(x_n)$, then
$$
(1+\g_n)\,h_n(w)<\rdm(1-|z|).
$$

\noindent{\rm (b)} If $\rlm(t)>\frac13\,\rlm(x_n)$, then
$$
2\bigl|\log|H_n'(w)|\bigr|+h_n(w)\le\frac32\,\rdm(1-|z|).
$$

\noindent{\rm (c)} If $\rlm(t)>\frac13\,\rlm(x_n)$ and
$|\arg(1-z)|>\frac12\pi/\rlm'(x_n)$, then
$$
2\log|H'_n(w)|+h_n(w)\le\frac{\rdm(1-|z|)}{1+\g_n}.
$$
\label{tl6}
\end{lemma}

\begin{proof} The first statement follows immediately from Lemma
\ref{tl5}. For convenience of notation, put
$$s=\log\frac1{1-|w|}.$$
We may rewrite estimate \f{41} as
\begin{equation}
\big|\exp[\rlm(t)]-\exp[\rlm(s)]\big|<1.
\label{47}
\end{equation}

\noindent(a) Using that $\rlm(t)\le\frac25\,\rlm(x_n)$ as well as the
convexity of $\rlm$, we obtain
$$\rlm(x_n)-(x_n-s)\rlm'(x_n)\le\rlm(s)\le\frac12\,\rlm(x_n),$$
and, as a consequence, $x_n-s\ge\frac12\,\rlm(x_n)/\rlm'(x_n)$.
Applying \f{100} and \f{101} with $z$ replaced by $w$, we complete the proof.
\smallskip

By the definition of $H_n$, using that $|1-w|\ge 1-|w|$ and
that, by the convexity of $\rlm$, $\rlm(x_n)+(s-x_n)\rlm'(x_n)\le
\rlm(s)$, we obtain:
\begin{equation}
\left.
\begin{gathered}
H_n'(w)=-\exp\big[\rlm(x_n)-x_n\rlm'(x_n)\big]\,\frac{F_{\rlm'(x_n)}(w)
\rlm'(x_n)}{1-w},\\
\log|H_n'(w)|\le\rlm(x_n)-x_n\rlm'(x_n)+(\rlm'(x_n)+1)\,s
+\log\rlm'(x_n),\\
\bigl|\log|H_n'(w)|\bigr|\le \rlm(s)+s+\log\rlm'(x_n).
\end{gathered}
\right\}
\label{44}
\end{equation}
In (b) and (c), we have $\rlm(t)>\frac13\,\rlm(x_n)$, and hence
$\rlm(s)>\frac14\,\rlm(x_n)$ because of \f{47}. By Lemma~\ref{tl1},
parts (b) and (c), we obtain
\begin{equation}
s+\log\rlm'(x_n)\le\left(2+\frac1{\e}\right)\log\rlm(s). 
\label{45}
\end{equation}
Moreover, as in \f{74}, we get
\begin{equation}
h_n(w)\le\exp[\rlm(s)].
\label{46}
\end{equation}
The assertion in part (b) now follows from \f{47}--\f{46}.
\smallskip

\noindent(c) We have
$$|\arg(1-z)|>\frac{\pi}{2\,\rlm'(x_n)}\quad\text{ and}\quad
|w-z|<e^{-2\rlm(t)}.$$
Since $\rlm(t)>\frac13\,\rlm(x_n)$, we get from Lemma~\ref{tl1}(c) that
$$e^{-2\rlm(t)}<\frac{e^{-2t}}{\rlm'(x_n)}.$$
A simple geometric argument then shows that
$$|\arg(1-w)|>\frac{\pi}{3\,\rlm'(x_n)}.$$
\smallskip

By Lemma~\ref{tl2}(d), and by the convexity of $\rlm$, we have
\begin{multline*}
\log(1+\g_n)+\log|h_n(w)|\le\rlm(x_n)+(s-x_n)\rlm'(x_n)+
\g_n-\frac{1}{3\rlm'(x_n)}\\
\le\rlm(s)+\g_n-\frac{1}{3\rlm'(x_n)}.
\end{multline*}
Applying once again Lemma~\ref{tl1}(c), and using \f{44} and \f{45},
we we are able to complete the proof.
\end{proof}

In the last technical lemma of this section, for nearby points $\xi$ and $z$
in the unit disk, we estimate the size of the quantities $|g_n'(\xi)/g_n(z)|$
and $|g_n(\xi)/g_n(z)|$, where $g_n$ is the zero-free analytic function
$$g_n(z)=\exp\left(\frac12\,H_n(z)\right).$$

\begin{lemma} In the notation of Lemma~{\rm\ref{tl4}}, take $\xi\in\D$
with $|\xi-z|<[\rdm(1-|z|)]^{-3}$. Then
\begin{gather*}
\left|\frac{g_n'(\xi)}{g_n(z)}\right|\le\exp\bigg[\frac{\rdm(1-|z|)}
{1+\g_n}-\rlm(t)\bigg],\\ 
\left|\frac{g_n(\xi)}{g_n(z)}\right|\le\exp
\bigg[\frac{\rdm(1-|z|)}{1+\g_n}-\rlm(t)\bigg].
\end{gather*}
\label{tl7}
\end{lemma}

\begin{proof} We prove only the first of these two inequalities; the
second one is treated analogously.

If $\rlm(t)\le\frac13\,\rlm(x_n)$, then by Lemma \ref{tl5},
for every $w$ with $|w-z|\le e^{-2\rlm(t)}$, we have
$\log\rdm(1-|w|)<\frac25\,\rlm(x_n)$, and by Lemma \ref{tl6},
$$
\sup_{|w-z|\le e^{-2\rlm(t)}}|g_n(w)|\le
\exp\bigg[\frac{\rdm(1-|z|)}{2(1+\g_n)}\bigg].
$$
Using the Cauchy integral formula, we then get
$$
|g_n'(\xi)|\le\exp\bigg[\frac{\rdm(1-|z|)}{2(1+\g_n)}+3\rlm(t)\bigg],
$$
By Lemma~\ref{tl4}(c), we obtain
\begin{equation*}
\left|\frac{g_n'(\xi)}{g_n(z)}\right|\le 
\exp\bigg[\frac{\rdm(1-|z|)}{1+\g_n}+3\rlm(t)-[\rlm(t)]^2\bigg]
\le\exp\biggl[\frac{\rdm(1-|z|)}{1+\g_n}-\rlm(t)\biggr]. 
\end{equation*}

If $\rlm(t)>\frac13\,\rlm(x_n)$, then we consider the following two
cases: $|g_n(z)|<1$ and $|g_n(z)|\ge1$.

If $|g_n(z)|<1$, then $|\arg(1-z)|\ge\frac12\pi/\rlm'(x_n)$, 
and by Lemma \ref{tl6}(c),
$$
|g_n'(\xi)|\le\exp\left[\frac{\rdm(1-|z|)}{2(1+\g_n)}\right].
$$
By Lemma~\ref{tl4}(c),
$$
|g_n(z)|\ge\exp\bigg[[\rlm(t)]^2-\frac{\rdm(1-|z|)}{2(1+\g_n)}\bigg],
$$
and we get 
$$
\left|\frac{g_n'(\xi)}{g_n(z)}\right|\le\exp
\bigg[\frac{\rdm(1-|z|)}{1+\g_n}-\rlm(t)\bigg]. 
$$

If, on the other hand, $1\le|g_n(z)|$, then by Lemma~\ref{tl6}(b),
$$
|g_n'(\xi)|\le\exp\bigg[\frac34\,\rdm(1-|z|)\bigg],
$$
so that
$$
\left|\frac{g_n'(\xi)}{g_n(z)}\right|\le\exp\bigg[\frac34\,\rdm(1-|z|)\bigg].
$$
This completes the proof of the lemma.
\end{proof}
\bigskip
\bigskip

\section{\label{main-constr}The construction of invertible noncyclic
functions}

Throughout this section, we suppose that $\bas$ satisfies \rf{10}
and $\bass$ is defined by \rf{11}. Using the technical results of Sections
\ref{aux-convex}--\ref{aux-harm-est}, in conjunction with Section
\ref{phrlind-est}, we produce here non-cyclic functions $F\in B^1(\D,\bas)$
satisfying various additional properties.

\begin{thm} There exists a zero-free function $F\in B^1(\D,\bas)$
such that $1/F$ is in
$B^1(\D,\bass)$ and $F^{1/p}$ is non-cyclic in $B^p(\D,\bas)$
for each $p$, $0<p<+\infty$.
\label{p1}
\end{thm}

\begin{thm} There exist functions $F_j$, $j=1,2,3,\ldots$, in
$B^1(\D,\bas)$, satisfying the conditions of
Theorem~{\rm\ref{p1}}, and such that for every $p$ with
$0<p<+\infty$ and for every integer $d$ with $1\le d\le+\infty$, the subspace
$[F_1^{1/p},\ldots,F_d^{1/p}]$ has index $d$ in $B^p(\D,\bas)$.
\label{p2}
\end{thm}

\begin{thm} There exists a function $F$ satisfying the
conditions of Theorem~{\rm\ref{p1}} and such that for
$f=F^{1/2}$ we have
$$
\int_\D\int_\D \frac1{|f(z)|^2}\,\frac{|f(z)-f(w)|^2}{|z-w|^2}\,
\bas(z)\bas(w)\mea(z)\mea(w)<+\infty.
$$
\label{p3}
\end{thm}

\begin{proofof}{Theorem \ref{p1}} We are going to
construct a harmonic function $V=\log|f|$ as the infinite sum of
functions harmonic in (different) neighborhoods of the
unit disk.

The setup is as in Section \ref{aux-convex} below. Starting with $\bas$,
we obtain
$$\bdm(s)=\log\frac1{\bas(1-s)}\quad\text{ and }\quad
\blm(x)=\log\bdm(e^{-x})$$
satisfying \rf{90} with some $\e$, $0<\e<1$. Then, applying
Lemma~\ref{tl1} to $\blm$, we obtain the minorant $\rlm$, and a sequence
$\{x_n\}_n$; we also have the function $\rdm$ defined by
$$\rdm(s)=\exp\left[\rlm\left(\log\frac1s\right)\right].$$
We fix $\kappa$ to be the quantity
$$\kappa=\exp\left(-1-\frac2\e\right).$$
We suppose that $r_n=1-e^{-x_n}$ tends
to $1$ as rapidly as required in Proposition~\ref{tl10}.

We set $V_1=0$ and argue by induction. At step $n$, we start with
a function $V_n$ harmonic in a neighborhood of the unit disk,
and find a real parameter $\eta$ with $0<\eta<1$, such that
\begin{gather}
|V_n(z)-V_n(w)|<1,\qquad |z-w|<\eta,\quad z,w\in\OD,
\label{0102}\\
\exp|V_n(z)|<\frac 1\eta,\qquad z\in\OD.
\label{31}
\end{gather}
Fix $\{x_n\}_n$ as a subsequence of the sequence in Lemma \ref{tl1}
which grows so fast that
\begin{equation}
e^{-x_n}<\eta\kappa\,e^{-2n},
\label{32}
\end{equation}
and put
$$\delta_n=e^{-x_n},\quad r_n=1-\delta_n,\quad\text{and}\quad
\g_n=e^{-\rlm(x_n)/10}=e^{-\blm(x_n)/10}.$$
We recall from Section \ref{phrlind-est} that the integer $N_n$ is such that
$$
N_n\le\frac{\kappa}{1-r_n}<N_n+1.
$$ 
By \f{31} and \f{32}, we get
$$
\delta_n^2\le\frac{\exp[-V_n(z)]}{n^2N_n}\le 1,\qquad z\in\OD.
$$
In Section \ref{aux-harm-constr}, the harmonic function $h_n$ is
constructed; it is given explicitly by formula (\ref{0101}). 
Using Lemma~\ref{tl4}(a) and Lemma~\ref{tl41}(c), for every $k$,
$0\le k<N_n$, we can choose numbers $\gamma_{n,k}$ with $0\le\gamma_{n,k}
\le\gamma_n$, such that 
\begin{equation}
\int_{\{z\in\D:\,|z-r_n|<\delta_n^2\}}\!\!\!\!
\exp\Big[\tilde h_{n,k}(z)-\bdm(1-|z|)\Big]\,
\mea(z)
=\frac{\exp\big[-V_n(e^{2\pi ik/N_n})\big]}{n^2N_n},
\label{33}
\end{equation}
where, as in Section~\ref{aux-harm-constr},
$$\tilde h_{n,k}(z)=(1+\gamma_{n,k})\,h_n(z).$$
We consider the following $N_n$ equidistributed points on the unit circle
$\T$, 
$$
\zeta_{n,k}=e^{2\pi ik/N_n},
$$ 
and as in Section~\ref{main-constr}, we put
$$
w_{n,k}=r_n\,\zeta_{n,k},\qquad 0\le k<N_n.
$$
Next, we consider the harmonic functions
\begin{align*}
U(z)&=\sum_{0\le k<N_n}h_{{n,k}}(z\bar\zeta_{n,k}),\\
V^0_{n+1}&=V_n+U.
\end{align*}
We need a little patch around each point $\zeta_{n,k}$ on the unit circle,
$$
\calO_{n,k}=\bigg\{z\in\D:|z-\zeta_{n,k}|<\delta_n\,e^{2/\e}\bigg\},
\qquad 0\le k<N_n.
$$
Note that these sets $\calO_{n,k}$, with $0\le k<N_n$, are mutually disjoint.
We also need the small disk around the point $w_{n,k}$ given by
$$\calD_{n,k}=\big\{z\in\D:\,|z-w_{n,k}|<\delta_n^2\big\},$$
where we recall that $\delta_n=1-r_n$. 

By Lemma~\ref{tl41}(a), for some constant $C$ depending only on
the parameter $\e$, we have 
\begin{equation}
\left.
\begin{gathered}
\big|U(z)-\tilde h_{{n,k}}\big(e^{-2\pi ik/N_n}z\big)\big|<C\,(1-r_n),
\qquad z\in \calO_{n,k},\\
|U(z)|<C\,(1-r_n),\qquad z\in\D\setminus \bigcup_k \calO_{n,k}.
\end{gathered}
\right\}
\label{106}
\end{equation}
In particular,
\begin{equation}
\!\!\big|V^0_{n+1}(z)-V_n(z)\big|=|U(z)|<C\,(1-r_n)\quad\text{ for}\quad
|z|\le1-\delta_n\,e^{2/\e}.
\label{35}
\end{equation}

Brought together, the relations \f{0102}, \f{33}, and \f{106} give us 
\begin{gather}
\int_{\calD_{n,k}}
\exp\Big[V^0_{n+1}(z)-\bdm(1-|z|)\Big]
\mea(z)\asymp\frac 1{n^2N_n}, \label{155}\\
\intertext{and by Lemma~\ref{tl4}(b), we obtain}
\int_{\{z\in\D:\,1-|z|<\delta_n\exp(2/\e)\}}
\exp\Big[V^0_{n+1}(z)-\bdm(1-|z|)\Big]\mea(z)
\asymp\frac 1{n^2},
\label{105}
\end{gather}
where we use the notation $a\asymp b$ for the relation
$a/c<b<c\,a$ with some positive constant $c$ depending only on $\e$.

The function $y\mapsto y-[\log y]^2$ increases monotonically as $y$ grows to
$+\infty$, whence we conclude that 
$$
\bdm(1-|z|)-\big[\log\bdm(1-|z|)\big]^2>
\rdm(1-|z|)-\big[\log\rdm(1-|z|)\big]^2,
$$  
and by Lemma~\ref{tl4}(c), \f{32}, and \f{106}, we have 
\begin{equation}
\int_{\{z\in\D:\,1-|z|<\delta_n\exp(2/\e)\}}
\exp\Big[-V^0_{n+1}(z)-\bdm(1-|z|)+\big[\log\bdm(1-|z|)\big]^2\Big]\mea(z)
<\frac 1{n^2}. 
\label{107}
\end{equation}
Moreover, Lemma~\ref{tl41}(b), for some $w\in\D$ with
$\delta_n/2<1-|w|<\delta_n$, we have
\begin{equation}
V^0_{n+1}(w)<-\frac23\,\rdm(\delta_n).
\label{108}
\end{equation}
Replacing $V^0_{n+1}$ by $V_{n+1}(z)=V_n(z)+U(\tau z)$ with
$\tau$ sufficiently close to $1$, $0<\tau<1$, we get the same
properties \f{35}--\f{108} with $V_{n+1}$ harmonic in a
neighborhood of the unit disk. 

As a consequence of \f{32} and \f{35}, the functions $V_n$
converge uniformly on compact subsets of the unit disk to a
harmonic function $V$ as $n\to+\infty$. We consider the
corresponding analytic function $F=\exp\big[V+i\widetilde{V}\big]$,
where the tilde indicates the harmonic conjugation operation,
normalized so that $\widetilde{V}(0)=0$.

It follows from \f{35} and \f{105} that
\begin{equation}
\int_\D|F(z)|\,e^{-\bdm(1-|z|)}\mea(z)<+\infty.
\label{a1}
\end{equation}
Analogously, by \f{35} and \f{107}, we have
\begin{equation}
\int_\D |F(z)|^{-1}e^{-\bdm(1-|z|)}
\exp\bigl([\log\bdm(1-|z|)]^2\bigr)
\mea(z)<+\infty,
\label{a2}
\end{equation}
and by \f{35} and \f{155}, we get for each $n=1,2,3,\ldots$,
\begin{equation}
\int_{\calD_{n,k}}|F(z)|\,e^{-\bdm(1-|z|)}\mea(z)\asymp\frac 1{n^2N_n},
\qquad k=0,1,2,\ldots,N_n-1. 
\label{36}
\end{equation}
By \f{35} and \f{108}, for every $n=1,2,3,\ldots$, there
exists a point $\xi_n\in\D $ with
$$\frac{\delta_n}2<1-|\xi_n|<\delta_n$$
such that
\begin{equation}
|F(\xi_n)|<\exp\left[-\frac12\,\rdm(\delta_n)\right].
\label{37}
\end{equation}

It remains to verify that $F^{1/p}$ is non-cyclic in $B^p(\bas)$ for each
$p$, $0<p<+\infty$. So, suppose that for a sequence of polynomials $q_j$,
$$
\big\|q_jF^{1/p}\big\|_{B^p(\bas)}^p=\int_\D|q_j(z)|^p\,|F(z)|\,
e^{-\bdm(1-|z|)}\,\mea(z)\le 1.
$$
The disks $\calD_{n,k}$ are disjoint for different indices $(n,k)$, so that
the above estimate has the immediate consequence
\begin{equation}
\sum_{n=1}^{+\infty}\sum_{k=0}^{N_n-1}\int_{\calD_{n,k}}
|q_j(z)|^p\,|F(z)|\,e^{-\bdm(1-|z|)}\,\mea(z)\le 1.
\label{sum-est}
\end{equation}
By the mean value theorem for integrals, there exist points
$z_{n,k}\in\calD_{n,k}$ (which may depend on the index $j$, too) such that
$$\int_{\calD_{n,k}}|q_j(z)|^p\,|F(z)|\,e^{-\bdm(1-|z|)}\mea(z)
=|q_j(z_{n,k})|^p\int_{\calD_{n,k}}|F(z)|\,e^{-\bdm(1-|z|)}\mea(z).$$
In view of \f{36}, we see that (\ref{sum-est}) leads to
$$\sum_{n=1}^{+\infty}\sum_{k=0}^{N_n-1}
\frac 1{n^2N_n}\,|q_j(z_{n,k})|^p\le C,$$
for some positive constant $C$, and by a weak type estimate, we obtain for
$n=1,2,3,\ldots$ that
$$
\sum_{k=0}^{N_n-1}|q_j(z_{n,k})|^p\le C\,n^2N_n.
$$
We find ourselves in the situation described in Section~\ref{phrlind-est}.
An application of Proposition~\ref{tl10} yields
$$
|q_j(z)|\le c\,\exp\left[\frac 1{1-|z|}\right],\qquad z\in\D.
$$
for some constant $c$ which is independent of the index $j$. The inequalities
\f{37} show that $q_jF^{1/p}$ cannot converge to a non-zero constant
uniformly on compact subsets of $\D$. As a result, $F^{1/p}$ is not
cyclic in $B^p(\bas)$.
\end{proofof}

\begin{proofof}{Theorem \ref{p2}} For the sake of
simplicity, we consider only the case $p=1$, $d=2$.
Arguing as in the proof of Theorem~\ref{p1}, we can produce 
invertible non-cyclic elements $F_1$ and $F_2$ in $B^1(\bas)$,
numbers
$$
\ldots<r_{1,n}<r_{2,n}<r_{1,n+1}<r_{2,n+1}<\ldots\to 1,\qquad 
n\to+\infty,
$$
integers $N_{j,n}$ such that
$$
N_{j,n}\le \frac\kappa{1-r_{j,n}}<N_{j,n+1},\qquad j=1,2,
$$ 
points
$$
w_{j,n,k}=r_{j,n}e^{2\pi ik/N_{j,n}},\qquad 0\le k<N_{j,n},
$$
and disks 
$$
\calD_{j,n,k}=\Big\{z\in\D:\,|z-w_{j,n,k}|<(1-r_{j,n})^2\Big\},\qquad
k=0,1,2,\ldots,N_{j,n}-1,
$$
with the properties that
\begin{equation}
\int_{\calD_{j,n,k}}|F_j(z)|\,e^{-\bdm(1-|z|)}\mea(z)
\asymp\frac 1{n^2N_{j,n}},\qquad j=1,2,
\label{42a}
\end{equation}
and
\begin{equation}
|F_{j}(z)|e^{-\rdm(1-|z|)}\le\exp\bigg[-\frac1{1-|z|}\bigg],
\quad z\in \calD_{3-j,n,k},\,\,\,j=1,2.\label{42b}
\end{equation}

If, for some polynomials $q_1$ and $q_2$, we are given that
\begin{equation}
\big\|q_1F_1+q_2F_2\big\|_{B^1(\bas)}\le 1,\label{42c}
\end{equation}
and
\begin{equation}
\|q_1\|_{H^\infty(\D)}+\|q_2\|_{H^\infty(\D)}\le A,\label{42d}
\end{equation}
then we argue as in the proof of Theorem~\rm\ref{p1}.
First, for sufficiently big $n$, $n_0=n_0(A)\le n<+\infty$, the
inequalities \f{42b}--\f{42d} imply that for $j=1,2$,
\begin{gather*}
\int_{\calD_{3-j,n,k}}\big|q_{j}(z)\,F_{j}(z)\big|\,e^{-\bdm(1-|z|)}
\mea(z)\le\frac 1{N_{3-j,n}},\\
\sum_{0\le k<N_{j,n}}\int_{\calD_{j,n,k}}\big|q_j(z)F_j(z)\big|
\,e^{-\bdm(1-|z|)}\mea(z)\le 2,
\end{gather*}
and hence, by \f{42a}, 
for some points $z_{j,n,k}$ in the disks $\calD_{j,n,k}$, we have
\begin{equation}
\sum_{0\le k<N_{j,n}}|q_j(z_{j,n,k})|\le 2n^2N_{j,n},\qquad n=n_0,n_0+1,
n_0+2,\ldots.
\label{42e}
\end{equation}
We have almost the assumption of Proposition~\rm\ref{tl10}, only with 
$p=1$ and the estimates start from index $n=n_0$. If we analyze the proof
of that proposition carefully, we obtain the corresponding estimate of
the function $q_j$,
\begin{equation}
|q_j(z)|\le c\,n_0^4\exp\left[\frac1{1-|z|}\right],\qquad z\in\D,
\label{42f}
\end{equation}
where $c=c(\kappa)$ is a positive constant, independent of the value of $n_0$. 
Similarly, the inequalities \f{42b}, \f{42c}, and \f{42f}
imply that the estimates $(j=1,2)$
\begin{gather*}
\int_{\calD_{3-j,n,k}}\big|q_{j}(z)\,F_{j}(z)\big|\,e^{-\bdm(1-|z|)}
\mea(z)\le\frac{c\,n_0^4}{R^2_{3-j,n}}\le\frac 1{N_{3-j,n}},\\
\sum_{0\le k<N_{j,n}}\int_{\calD_{j,n,k}}\big|q_j(z)\,F_j(z)\big|
\,e^{-\bdm(1-|z|)}\mea(z)\le 2,
\end{gather*}
hold for all $n$ such that
\begin{equation}
c\,n_0^4\le \min\big\{N_{1,n},N_{2,n}\big\}.
\label{42m}
\end{equation}
As a consequence of \f{32}, $N_{j,n}\ge e^n$, and the inequality
\f{42m} holds for all $n\ge n_1$, with $n_1$ equal to the integer part of
$\alpha\log n_0$, for some positive real parameter $\alpha$; clearly, for
sufficiently big $n_0$, we have $n_1<n_0$.
Arguing as before, we get \f{42e} for all $n\ge n_1$, and then
$$
|q_j(z)|\le c\,n_1^4\exp\left(\frac 1{1-|z|}\right),\qquad z\in\D.
$$
Continuing in this way, we get \f{42e} for all 
$n\ge n_\infty=\lim_{k\to\infty}n_k$, with $n_\infty$ 
independent of $A$, and, as a consequence,
$$
|q_1(z)|+|q_2(z)|<c\,\exp\left(\frac 1{1-|z|}\right),\qquad z\in\D,
$$
for some constant $c$, which is independent of $q_1$, $q_2$, and $A$.
This shows that (see, for instance, \cite{R} for $d=2$, and 
\cite[Section~3]{Bor} for the general case) 
$$
\ind([F_1,F_2])=2.  
$$
The proof is complete.
\end{proofof}

\begin{proofof}{Theorem \ref{p3}} Let $F$ be the
function constructed in the proof of Theorem~\ref{p1}, and
consider $f=F^{1/2}$. The inequalities \f{a1} and \f{a2} together
with the inequality
$$
x^6<\exp\bigl[(\log x)^2\bigr]\quad\text{ for }\,\, x>e^6,
$$
and the identity
$$\bas(z)=e^{-\bdm(1-|z|)}$$
show that we have 
\begin{gather*}
\int_\D|f(z)|^2\bas(z)\mea(z)<+\infty,\\
\int_\D \frac1{|f(z)|^2}[\bdm(1-|z|)]^6\bas(z)\mea(z)<+\infty.
\end{gather*}

Let $\calE$ be the set
$$\calE=\Big\{(z,w)\in\D\times\D:\,|w-z|\le[\rdm(1-|z|)]^{-3}\Big\},$$
which is a rather small neighborhood of the diagonal; we then have the
estimate
\begin{multline*}
\iint_{\D^2\setminus\calE}
\,\,\frac1{|f(z)|^2}\frac{|f(z)-f(w)|^2}{|z-w|^2}\,
\bas(z)\,\bas(w)\mea(z)\mea(w)\\
\le 2\iint_{\D^2}\frac{|f(z)|^2+|f(w)|^2}{|f(z)|^2}\,
\big[\rdm(1-|z|)\big]^6 \bas(z)\,\bas(w)\mea(z)\mea(w)\\
\le C\iint_{\D^2}\frac{1+|f(z)|^2}{|f(z)|^2}\,\big[\rdm(1-|z|)\big]^6
\bas(z)\mea(z)<+\infty.
\end{multline*}
By Lemma \ref{tl6},
\begin{equation}
\big|\rdm(1-|z|)-\rdm(1-|w|)\big|<1\,\,\hbox{ for }\,\,
(z,w)\in\calE.
\label{est-000}
\end{equation}
Hence, it suffices to verify that for $(z,w)\in\calE$, we have 
\begin{equation}
\frac1{|f(z)|^2}\frac{|f(z)-f(w)|^2}{|z-w|^2}\le M
\exp\big[2\rdm(1-|z|)\big],
\label{est-001}
\end{equation}
for some positive constant $M$. This is so because, as a consequence of
(\ref{est-001}), we have the estimate
\begin{multline*}
\iint_{\calE}
\,\,\frac1{|f(z)|^2}\frac{|f(z)-f(w)|^2}{|z-w|^2}\,
\bas(z)\,\bas(w)\mea(z)\mea(w)\\
\le M\iint_{\calE}\exp\big[2\rdm(1-|z|)\big]\,
\bas(z)\,\bas(w)\mea(z)\mea(w)\\
=M\iint_{\calE}
\exp\Big[2\rdm(1-|z|)-\bdm(1-|z|)-\bdm(1-|w|)\Big]\,
\mea(z)\mea(w)\\
\le M\iint_{\calE}\exp\Big[\rdm(1-|z|)-\rdm(1-|w|)\Big]\,
\mea(z)\mea(w),
\end{multline*}
and the latter is bounded by $M$, if we use (\ref{est-000}). We turn to the
verification of (\ref{est-001}).
We find that there exists a point $\xi=\xi(z,w)\in\D$ with
$$|\xi-z|\le[\rdm(1-|z|)]^{-3},$$
such that
$$
\left|\frac{f(z)-f(w)}{z-w}\right|\le |f'(\xi)|,
$$
by expressing the function $f(z)-f(w)$ as a path integral from $z$ to $w$.
It follows that we need only to verify that
$$
\frac{|f'(\xi)|^2}{|f(z)|^2}\le M\,\exp\big[2\rdm(1-|z|)\big]. 
$$
Now, we fix $\zeta$ and $\xi$ in $\D$, with
$$|\zeta-\xi|\le[\rdm(1-|\zeta|)]^{-3},$$
so that we are in the setting above (only with slightly different variable
names). We recall the details of the construction of the function $F$ in
the proof of Theorem~\ref{p1}. We then find a positive integer $n=n(|\zeta|)$
such that
\begin{equation}
\delta_{n+1}\,e^{2/\e}=\exp\left(-x_{n+1}+\frac2\e\right)\le1-|\zeta|
<\exp\left(-x_{n}+\frac2\e\right)=\delta_{n}\,e^{2/\e};
\label{59a}
\end{equation}
this is definitely possible at least if $\zeta$ is reasonably close to $\T$.
On the other hand, if we start with a positive integer $n$, we may form
the annulus
\begin{equation}
\calU_n=\Big\{z\in\D:\,\delta_{n+1}\,e^{2/\e}<1-|z|\le
\delta_{n}\,e^{2/\e}\Big\},
\label{59a'}
\end{equation}
so the above relation (\ref{59a}) places the point $\zeta$ inside $\calU_n$.
Let $V_n$ be the harmonic function appearing in the proof of Theorem \ref{p1},
and let $W_n=V_n+i\widetilde{V}_n$ be the analytic function having $V_n$ as
real part. We form the exponentiated function
$$F_n(z)=\exp\big(W_n(z)\big),\qquad z\in\D,$$
which is analytic and zero-free in a neighborhood of the closed disk $\OD$. 
The construction of the function $V_n$ (and hence that of $W_n$) involves
the choice of the points $x_1,x_2,\ldots,x_{n-1}$.
By letting the sequence $\{x_n\}_n$ tend to $+\infty$ as rapidly
as need be, we can make sure that 
\begin{equation}
\max_{z\in\OD}\left\{|F_n(z)|,\frac1{|F_n(z)|},
|F_n^\prime(z)|\right\}\le\rdm\big(e^{-x_n+2/\e}\big)\le
\exp\Big(\frac14\,\lambda(t)\Big),
\label{59b}
\end{equation}
where
$$t=\log\frac1{1-|\zeta|},$$
and the right-hand inequality in (\ref{59b}) holds because of (\ref{59a}).
The function $F$ is the limit of the functions $F_n$ as $n\to+\infty$,
and it is of interest to understand the ``tail'' function $F/F_n$. 
By Lemma~\ref{tl41}(a) and some elementary estimates of harmonic functions
and their gradients, plus the fact that the points $x_n$ approach $+\infty$
very rapidly as $n\to+\infty$, the inequalities \f{59a} imply that 
\begin{equation}
\max_{z\in\overline\calU_n}\left\{\bigg|\frac{F}{F_{n+1}}(z)\bigg|,
\bigg|\frac{F_{n+1}}{F}(z)\bigg|,\bigg|\left(\frac{F}{F_{n+1}}\right)'(z)
\bigg|\right\}\le C, \label{59c}
\end{equation}
for some positive constant $C$ independent of $n$; the notation
$\overline\calU_n$ stands for the closure of $\calU_n$. We need to apply
the estimate (\ref{59c}) to the two points $\zeta$ and $\xi$, and although
$\xi$, strictly speaking, need not belong to $\calU_n$, it is not far away
from this set, and we can make sure that the estimate (\ref{59c}) holds for
it, simply because Lemma ~\ref{tl41}(a) applies in a slightly bigger annulus
than $\calU_n$.
Next, we pick the point $\zeta_{n,k_0}$, with $0\le k<N_{n}$, which
is closest to the given point $\zeta\in\D$; after a rotation of the disk,
we may assume $\zeta_{n,k_0}=\zeta_{n,0}=1$. 
Since $|\zeta-\xi|$ is much smaller than $1-|\zeta|$, we obtain
$$
\min_{k\ne 0}\big|\zeta_{n,k}-\zeta\big|
\ge\exp\left[-x_n+\frac2\e\right].
$$
Let $G_n$ be the function
$$G_n(z)=\frac{F_{n+1}(z)}{F_n(z)}\,\exp\Big[-(1+\gamma_{n,0})\,H_n(z)\Big],$$
where $\gamma_{n,0}$ is determined by \f{33}, and the analytic function
$H_n$ is defined in \f{0101}; the real part of $H_n$ equals $h_n$.
Then, again by Lemma~\ref{tl41}(a),
\begin{equation}
\max_{z\in\{\zeta,\xi\}}\max\Big\{|G_n(z)|,\frac1{|G_n(z)|},
|G^\prime_n(z)|\Big\}\le C, \label{59d}
\end{equation}
where $C$ is some positive constant which does not depend on $n$. 
In conclusion, we write
$$f(z)^2=F(z)=g_n(z)^2\,F_n(z)\,\frac{F}{F_{n+1}}(z)\,G_n(z),$$
where
$$g_n(z)=\exp\left(\frac12\,\big(1+\gamma_{n,0}\big)\,H_n(z)\right),$$
and by Lemma~\ref{tl7} (note that the function $g_n$ appearing in the
lemma is slightly different),
\begin{equation}
\left.
\begin{aligned}
\frac{|g_n'(\xi)|^2}{|g_n(\zeta)|^2}&\le
\exp\bigl[2\rdm(1-|\zeta|)-2\rlm(t)\bigr],\\
\frac{|g_n(\xi)|^2}{|g_n(\zeta)|^2}&\le
\exp\bigl[2\rdm(1-|\zeta|)-2\rlm(t)\bigr].
\end{aligned}
\right\}
\label{59e}
\end{equation}
Hence, in view of \f{59b}--\f{59e}, we obtain
$$
\frac{|f'(\xi)|^2}{|f(\zeta)|^2}\le Ce^{\lambda(t)}
\frac{|g_n(\xi)|^2+|g_n'(\xi)|^2}{|g_n(\zeta)|^2}
\le\exp\big[2\rdm(1-|\zeta|)\big]. 
$$
The proof is complete.
\end{proofof}

\bigskip

\noindent \textsc{Alexander Borichev, Department of Mathematics,
University
of Bordeaux I, 351, cours de la Lib\'eration, 33405 Talence, France}

\noindent\textsl{E-mail}: \texttt{borichev@math.u-bordeaux.fr}
\medskip

\noindent \textsc{H\aa{}kan Hedenmalm, Department of Mathematics,
The Royal Institute of Technology, S--100 44 Stockholm, Sweden}

\noindent\textsl{E-mail}: \texttt{haakanh@math.kth.se}
\medskip

\noindent \textsc{Alexander Volberg, Department of Mathematics, 
Michigan State University, East Lansing, MI 48824, USA}

\noindent\textsl{E-mail}: \texttt{volberg@math.msu.edu}
\medskip
\enddocument
\bye